\documentclass[12pt]{iopart}
\usepackage[english]{babel}
\expandafter\let\csname equation*\endcsname\relax
\expandafter\let\csname endequation*\endcsname\relax
\usepackage{amsmath}
\usepackage{amsfonts}
\usepackage{bm}
\DeclareMathOperator{\cn}{cn}
\DeclareMathOperator{\sn}{sn}
\DeclareMathOperator{\dn}{dn}
\DeclareMathOperator{\eK}{K}

\DeclareMathOperator*{\res}{Res}
\newcommand\Ord{\Or}
\newcommand\R{\hat{\boldsymbol{\mathcal R}}}
\newcommand\Sr{\hat{\boldsymbol{\mathsf R}}}
\newcommand\Lo{\hat{\boldsymbol{\mathsf L}}}
\newcommand\Po{\hat{\boldsymbol{\mathsf  P}}}
\newcommand\So{\hat{\boldsymbol{\mathsf  S}}}
\newcommand\Us{\hat{\boldsymbol{\mathsf  U}}}
\newcommand\Vs{\hat{\boldsymbol{\mathsf  V}}}
\newcommand\Z{\hat{\boldsymbol{\mathsf  Z}}}

\newcommand\D{\hat{\boldsymbol{\mathsf  D}}}
\newcommand\expD[1]{{\rm exp}_{\hskip-0.1em D}(#1)}
\newcommand\iexpD[1]{{\rm exp}_{\hskip-0.1em D}^{-1}(#1)}
\newcommand\vecp[1]{\vec{#1}}
\newcommand\Checked{}

\begin{document}
\title[Kato perturbation expansion in classical mechanics]{Kato perturbation expansion in classical mechanics and an explicit expression for a Deprit generator.}

\author{Nikolaev~A.\,S.$^{1,2}$}
\address{$^1$  Institute of Computing for Physics and Technology, 6, Zavodskoy proezd, Protvino, Moscow reg., Russia, 142281}
\address{$^2$  RDTeX LTD, 1 bld. 40, Nagatinskaya Str., Moscow, Russia, 117105}
\ead{Andrey.Nikolaev@rdtex.ru, http://andreynikolaev.wordpress.com}

\begin{abstract}
{This work explores the structure of Poincare-Lindstedt perturbation series in Deprit operator formalism and establishes its connection to Kato resolvent expansion.
A discussion of invariant definitions for averaging and integrating perturbation operators and their canonical identities reveals a regular pattern in a Deprit generator. 
The pattern was explained using  Kato series and the relation of   perturbation  operators to Laurent coefficients for the resolvent of Liouville operator.

This purely canonical approach systematizes the series and leads to the  explicit expression  for the Deprit generator in any perturbation order:
\[G = - \hat{\boldsymbol{\mathsf  S}}_H H_i.\]
Here, $\hat{\boldsymbol{\mathsf  S}}_H$  is the  partial pseudo-inverse of the  perturbed Liouville operator.  
Corresponding Kato series provides a reasonably effective computational algorithm. 

The canonical connection of perturbed and unperturbed averaging operators allows for a description of  
  ambiguities in the generator and transformed Hamiltonian, while Gustavson integrals turn out to be   insensitive to normalization style. 
Non-perturbative examples  are used for illustration.
}
\end{abstract}


\section{Introduction.}
This work was inspired by remarkable analogies between mathematical formalisms of perturbation expansions in classical and quantum mechanics. For instance,
 classical secular perturbation theory  \cite{Cary} corresponds  to time-dependent quantum mechanical perturbation expansion \cite{messiah}. 
``Action-angle'' variables correspond to energy representation in quantum mechanics.
The classical Poincare-Lindstedt method \cite{Poincare} in  Lie algebraic formalism \cite{Deprit,DragtFinn} has a direct analogy in quantum Van Vleck perturbation expansion  \cite{VanVleck,Shavitt,Klein}. Similarly, classical Birkhoff-Gustavson normal forms relate to 
quantum mechanical perturbation theory in Bargmann-Fock space \cite{Graffi2,Nikolaev}.

However, quantum mechanics has a wider diversity of perturbation methods. Some of these methods may be of interest for classical perturbation theory.
 Here we will construct a classical analogue of Kato series \cite{Kato}.

Kato used the Laurent and Neumann expansions of a resolvent  operator around the eigenvalues of the quantum mechanical Hamiltonian.
In classical mechanics, our tools will be the Liouville operator and its resolvent.
But the spectrum of classical Liouvillian  is more complex than  
that of quantum Hamiltonian \cite{spohn75}. 
Also, the multidimensional classical perturbation series must diverge on everywhere dense set of resonances \cite{Arnold}. 

Despite  this internal divergence, canonical perturbation expansion is an efficient tool
in celestial mechanics  \cite{FM},  nonlinear physics and accelerator theory.
This is why new perturbative algorithms are  important.
We  focus here on formal constructions of canonical perturbation series  and general formulae. 
Their convergence and nonresonance cancellations will be discussed in another article. 

It is worth noting that noncanonical perturbation expansions of Liouvillian resolvents were  used by 
the Brussels-Austin Group in nonequlibrium statistical mechanics  \cite{Balescu75}. However,  we  will construct purely canonical  series.

This paper is organized as follows. First, we  review the construction of 
Poincare-Lindstedt-Deprit series using  invariant operator formalism and introduce
convenient notation borrowed from quantum mechanics. 
 We then  discuss algebraic properties of basic perturbation operators. 
 A set of canonical identities  uncovers regular pattern in Deprit generator.

To extend  this pattern to all orders,  we   borrow { Abel averaging}  from quantum mechanics and establish a connection 
to the resolvent of Liouville operator and Kato  expansion. 
Canonical properties of Liouvillian resolvent   allow us to write a simple explicit expression for Deprit generator: 
\begin{displaymath}
G = - \So_H H_i.
\end{displaymath}
In this formula, the integrating operator  $\So_H$  is the partial pseudo-inverse of the {\em perturbed} Liouville operator.

After providing a description of the generator structure in any perturbation order  and its ambiguity and difference in normalization style from the standard Deprit algorithm, 
we will extend  the formulae to a multidimensional case and discuss Gustavson integrals and computational efficiency. Finally,  we will illustrate the findings using non-perturbative examples.

Our method combines classical and quantum mechanical perturbation approaches. The novel results are the observation of a regular pattern in the perturbation series for the generator, its explanation using Kato series, the explicit expression for Deprit generator in all perturbation orders and  the  insensitivity of Gustavson integrals to normalization style. 

The 
``Supplementary data files'' contain demonstrations  and large formulae, including the general normal form of Hamiltonian  up to  
the seventh perturbation order. The demonstrations use the freeware computer algebra system FORM \cite{Form}.

\section{Classical perturbation expansion.}
\subsection{Liouville operator.}

Equations of motion, for most dynamical systems, cannot be solved analytically. Therefore, we should be interested in approximations that  preserve underlying physical structures. One   such approximation is the classical perturbation theory \cite{Arnold}.
It explores the dynamical behaviour of {\bf d}-dimensional mechanical system with Hamiltonian, which differs from
the  solvable (integrable) system by the small perturbation
\[
H=H_0+\alpha H_i,
\]
where $H_0$ and $H_i$ are the functions of canonical variables
$\vecp p$, $\vecp q$ on $\mathbb{R}^d\!\times\!\mathbb{R}^d$.
We will consider only autonomous (time-independent) Hamiltonians having compact energy surfaces $H(\vecp p,\vecp q)=E$.
We also  assume that all  functions  are analytic. 

Perturbation theory approximates the time evolution of canonical variables 
and constructs  integrals of motion for perturbed system. More generally, it approximates the time evolution of function $F(\vecp p,\vecp q)$, obeying Hamiltonian equations:
\[
\frac{\rmd F}{\rmd t} = [F,H],\qquad [F,H] = \sum_{i=1}^{\rm\bf d}
\frac{\partial F}{\partial q_i}\frac{\partial H}{\partial p_i} -
\frac{\partial F}{\partial p_i}\frac{\partial H}{\partial q_i}.
\]
Such evolution defines  the continuous one-parametric family of canonical  transformations of the phase space  named {\it Hamiltonian flow}. Using  Liouville operator
\[
\Lo_{H} = [\,.\,,H], \qquad \Lo_{H}F = [F,H],
\]
the formal solution of Hamiltonian equations 
for an autonomous Hamiltonian may be written as operator exponent
\begin{equation}
{\left. \widetilde F( p, q)\right|}_t 
= {\left. \rme^{t \,\Lo_{H}} F(p,q)\right|}_0,\quad \Lo_{H} \rme^{t \,\Lo_{H}}= \rme^{t \,\Lo_{H}}\Lo_{H}.\label{u4}
\end{equation}
This is  a canonical  transformation commutative with the autonomous $\Lo_{H}$. 

\subsection{Near-identity canonical transformations.}
 
Classical perturbation theory uses a more general construction of a family of 
near-identity $\alpha$ dependent canonical transformations on  
$\mathbb{R}^{2{\bf d}}$ phase space 
\[
\widetilde {\bf x} = \Us(\alpha)  {\bf x},
\qquad \begin{cases} x_i&=q_i,\\ x_{i+{\rm\bf d}}&=p_i,\end{cases}\quad i=1,\,\ldots\,
{\rm\bf d},
\]
equipped with the canonical structure
\[
[F({\bf x}),G({\bf x})] = \sum_{i,j}^{\rm\bf2d}
\frac{\partial F}{\partial x_i} \omega_{ij} \frac{\partial G}{\partial x_j}. 
\]
Here, $\omega_{ij}= - \omega_{ji}$ is the non-degenerate skew-symmetric Jacobi matrix.
Such canonical transformations  are always Hamiltonian flows \cite{Deprit,Giorgilli:arXiv1211.5674} 
\[\frac{\partial\Us}{\partial\alpha} = \Lo_G \Us\]
 in ``time'' $\alpha$ \ with some generator $G({\bf x},\alpha)$. The perturbation theory constructs this generator as the power series  
$G({\bf x},\alpha) = \sum_0^\infty \alpha^n G_n({\bf x}).$
It also computes the  series for the transformation $\Us=\sum_0^\infty \alpha^n \Us_n$ and its inverse, $\Vs=\Us^{-1}$, which satisfies the equation
\begin{equation}
\frac{\partial\Vs}{\partial\alpha} = - \Vs \Lo_G.\label{inverseD}
\end{equation}
In  celestial mechanics these formulae 
are  known as   ``Lie transforms'' by Deprit \cite{Deprit}. 
Substituting  series for $G$, $\Us$ and $\Vs$ into  the above  equations, 
Deprit  obtained the recursive relations for coefficients: 
\[
\Us_n =  \frac{1}{n} \sum_{k=0}^{n-1} \Lo_{G_{n-k-1}} \Us_k, \qquad
\Vs_n = - \frac{1}{n} \sum_{k=0}^{n-1} \Vs_k \Lo_{G_{n-k-1}}. 
\]
These expressions can be iterated to derive the non-recursive formulae \cite{Cary,Koseleff} 
\begin{align*}
\Us_n&= \sum_{\substack{(m_1,\ldots,m_r)\\n>m_1>m_2>\cdots>m_r}} \;\;\frac{\Lo_{G_{n-m_1-1}}}{n} \frac{\Lo_{G_{m_2-m_1-1}}}{m_1}\cdots\frac{\Lo_{G_{m_r-1}}}{m_r}\ ,\\
\Vs_n&= \sum_{\substack{(m_1,\ldots,m_r)\\n>m_1>m_2>\cdots>m_r}} 
\!\!\!\!\!\!(-1)^{r+1}\frac{\Lo_{G_{m_r-1}}}{m_r} \cdots\frac{\Lo_{G_{m_2-m_1-1}}}{m_1}\frac{\Lo_{G_{n-m_1-1}}}{n}\ .
\end{align*}\Checked
The sum runs  over all sets of 
integers $(m_1,\ldots,m_r)$, satisfying $n>m_1>\cdots>m_r>0$.

Hereafter, we will use the notation $\Us=\expD{\alpha\Lo_G}$,  $\Vs=\iexpD{\alpha\Lo_G}$ and term  ``Deprit exponents'' for these series to emphasize their exponential-like role.
Look at the first few orders:
\begin{multline}
\expD{\alpha\Lo_G} = 1 + \alpha \Lo_{G_0} + \frac{\alpha^2}{2}(\Lo_{G_0}^2 + \Lo_{G_1})\\
\shoveright{+ \frac{\alpha^3}{6} ( \Lo_{G_0}^3 + \Lo_{G_0}\Lo_{G_1} + 2 \Lo_{G_1}\Lo_{G_0} + 2 \Lo_{G_2}) 
+  \Ord(\alpha^3)}\\
\shoveleft{\iexpD{\alpha\Lo_G} = 1 -\alpha \Lo_{G_0} + \frac{\alpha^2}{2}(\Lo_{G_0}^2 -\Lo_{G_1})} \\
- \frac{\alpha^3}{6} ( \Lo_{G_0}^3 -
 2\Lo_{G_0}\Lo_{G_1} -  \Lo_{G_1}\Lo_{G_0} 
+ 2 \Lo_{G_2}) +  \Ord(\alpha^3)  
\end{multline}\Checked
See also   ``{\tt deprit\_exponents}'' demo in the supplementary data files. 


We will outline here  the key properties of canonical transformation $\Us$. This is the point transformation of phase space $\Us F({\bf x})=F( \Us {\bf x})$,
which also preserves  Poisson brackets.  It make sense to consider the Poisson bracket as an additional antisymmetric  product ({\em Poisson algebra}). Canonical transformations preserve both the algebraic and canonical structures of function space
\begin{align*}
\Us (FH)\ \ \ &=(\Us F)(\Us H),\\
\Us \left([F,H]\right)&=[\Us F,\Us H],
\end{align*}
for any functions  $F({\bf x})$ and $H({\bf x})$. 

Due to  the Leibniz rule and Jacobi identity, the Liouville operator is the ``derivation'' of both products 
\begin{align*}
\Lo_G(FH)&=(\Lo_G F)H+F(\Lo_G H),\\
\Lo_G [F,H]&= [\Lo_G F,H]+[ F,\Lo_G H].
\end{align*}\Checked
This dialgebraic property will be used to establish identities for perturbation operators. 

The operator formalism of generators, or {\em Lie-algebraic approach}, was introduced in the perturbation theory of classical  mechanics in 1966 \cite{Hori}.
This formalism is much simpler for computer algebra calculations than the previous ones that use {\em generating functions}.
It is interesting that  corresponding quantum mechanical formalism \cite{VanVleck}  appeared more than 40 years before its classical analogue.

Our goal is to transform canonically the perturbed Hamiltonian into a simpler,
preferably integrable form. 
The transformed Hamiltonian will become  an integral of motion of the unperturbed system.
Here we briefly review this standard construction. 

The canonical map   
$ \tilde{\bf{x}} = \expD{\alpha\Lo_G}\,{\bf{x}}$ transforms the perturbed Hamiltonian into
\begin{multline}
\widetilde H=  \iexpD{\alpha\Lo_G} (H_0+\alpha H_i) =H_0+\alpha (H_i-\Lo_{G_0}H_0)
+\Ord(\alpha^2)
\\=H_0+\alpha(\Lo_{H_0} G_0 + H_i)
+\frac{\alpha^2}{2} (\Lo_{H_0} G_1 + \Lo_{G_0}^2 H_0 - 2\Lo_{G_0} H_i  )
+\Ord(\alpha^3).\label{Htrans}
\end{multline}\Checked
To simplify coefficients, one needs to  
invert the Liouville operator $\Lo_{H_0}$.
Since this operator has a non-empty kernel, only  its pseudo-inverse can be constructed. 

\subsection{Basic perturbation operators.}

We assume that the unperturbed Hamiltonian system with the compact energy surface $H_0({\bf{x}})=E$ is completely integrable in a Liouville sense and performs   quasi-periodic motion on the invariant tori  \cite{Arnold}.  The time average 
\begin{equation}
\langle F\rangle =
\lim_{T\to\infty}\frac{1}{T} \int_0^{T} {\left. F(\vecp p(t),\vecp q(t))\right|}_{p(0)=p\atop q(0)=q} \,\rmd t = 
\lim_{T\to\infty}\frac{1}{T} \int_0^{T}
\rme^{t\Lo_{H_0}} F({\bf{x}}) \rmd t\label{Cesaro}
\end{equation}
exists in an ``Action--Angle'' representation 
for any analytic function $F(\bf{x})$. 
This average  is a function of the initial point ${\bf x}$.
Being written in an invariant form, it also exists in  other canonical variables.

The averaging  operation extracts from $F(\bf{x})$  its {\em secular} non-oscillating part, which remains constant under the time evolution. 
In functional analysis, it is known as  {\em Ces\`aro $(C,1)$   averaging} \cite{HP}. 
Corresponding  operator
\begin{equation}
\Po_{H_0} = \lim_{T\to\infty}\frac{1}{T} \int_0^{T} \rmd t\, \rme^{t\Lo_{H_0}}   \label{Pdef}
\end{equation}
 is projector $\Po_{H_0}^2 = \Po_{H_0}$. It projects $F({\bf x})$ onto the secular space of functions commutative with $H_0$. 
This is the kernel of $\Lo_{H_0}$ and the algebra of its integrals of motion.

The notation $\Po_{H_0}$ was borrowed from  quantum mechanics. Standard notations like ${\pmb{\langle\ \rangle}}$, $\bf\overline{F}$ and $M_t$ \cite{mitropolsky71,bogolyubov}  do not emphasize
the pattern in perturbation expansion that we  seek. 

Complementary projector $1-\Po_{H_0}$ extracts the time-oscillating part from $F({\bf x})$. 
It projects on the non-secular space of oscillating functions where the inverse of $\Lo_{H_0}$ exists (we will consider here only semi-simple $\Lo_{H_0}$).

This inverse is the integrating operator $\So_{H_0}$,  which is also known as the ``solution of homological equation'', ``tilde 
operation'', ``zero-mean antiderivative'' \cite{Sanders},  
``Friedrichs ${\bf \widehat\Gamma}$ operation'' \cite{Friedrichs}, ``${\bf 1/{k}}$ operator'', etc. Its invariant definition is \cite{mitropolsky71}
\begin{equation}
\So_{H_0} = - \lim_{T\to\infty} \frac{1}{T} \int_0^{T} \rmd t \int_0^t \rmd\tau \,\rme^{\tau \Lo_{H_0}} \left(1-\Po_{H_0}\right).\label{Sdef}
\end{equation}\Checked
See also \cite{cushman84,vorobiev11} for periodic extensions.  The formal calculation 
\begin{multline*} \Lo_{H_0} \So_{H_0} = - \lim_{T\to\infty} \frac{1}{T} \int_0^{T} \rmd t \int_0^t \rmd\tau \,
 \frac{\partial}{\partial \tau} \rme^{\tau\Lo_{H_0}} (1-\Po_{H_0}) 
\\= - \lim_{T\to\infty}\frac{1}{T} \int_0^{T} \rmd t\, (\rme^{t\Lo_{H_0}}-1) (1-\Po_{H_0}) = 1-\Po_{H_0}, 
\end{multline*}
confirms that $\So_{H_0}$ is the partial pseudo-inverse of the unperturbed Liouville operator:
\begin{align*}
\Lo_{H_0}\So_{H_0} &= \So_{H_0}\Lo_{H_0}=1-\Po_{H_0},\\
\So_{H_0}\Po_{H_0} &= \Po_{H_0}\So_{H_0}\equiv 0.
\end{align*}\Checked

These three basic operators $(\Lo_{H_i}$, $\Po_{H_0}$,  $\So_{H_0})$ are the building blocks of perturbation expansion.
Classical perturbation theory provides several algorithms to compute them:
\begin{itemize}
\item{\bf ``Action-Angle'' representation} \cite{Arnold}. With its roots deep in the $19\textsuperscript{th}$ century, 
 it can be traced to Le Verrier and Delaunay. In ``Action-Angle'' canonical coordinates, the unperturbed Hamiltonian 
is the function of  {\bf d} ``action'' variables $\vecp J$ only, and the perturbation is  $2\pi$ periodic in the phases $\vecp\phi$
\begin{displaymath}
H = H_0(\vecp J) + \alpha H_i(\vecp J,\vecp\phi) 
= H_0(\vecp J) + \alpha \sum_{\vecp k} \tilde H_i(\vecp J,\vecp k)\,  \rme^{\rmi(\vecp k,\vecp\phi)}.
\end{displaymath}
The motion of an unperturbed system is quasi-periodic: 
\[\vecp J = \textrm{const},\qquad\vecp\phi(t) = \vecp\omega t+\vecp\phi_0,\qquad\vecp\omega=\frac{\partial H_0}{\partial\vecp J}.
\] Because the Fourier components in the ``Action-Angle'' representation are eigenfunctions of 
$\Lo_{H_0}$, 
 the perturbation operators can be  written as
\begin{align}
\Po_{H_0} F &= \lim_{T\to\infty} \frac{1}{T} \int_0^{T} \rmd t \,
\sum_{\vecp k} \tilde F(\vecp J,\vecp k) \,
\rme^{\rmi(\vecp k,\vecp\phi_0)+ \rmi(\vecp\omega,\vecp k)t}\nonumber\\&
= \sum_{(\vecp\omega,\vecp k)=0} \tilde F(\vecp J,\vecp k) \,
\rme^{\rmi(\vecp k,\vecp\phi_0)},\label{AAdefs}\\
\So_{H_0} F &= \sum_{(\vecp\omega,\vecp k)\ne 0} \frac{1}{\rmi(\vecp\omega,\vecp k)}
\tilde F(\vecp J,\vecp k) \, \rme^{\rmi(\vecp k,\vecp\phi_0)}.\nonumber
\end{align}
These  well-known expressions are frequently used as the definitions of $\Po_{H_0}$ and $\So_{H_0}$.

\item{\bf ``Birkhoff-Gustavson-Bruno'' normalization}  \cite{Birkhoff,Gustavson,Bruno} for power series. In its simplest case, a quadratic 
unperturbed Hamiltonian is diagonalizable into 
\[H_0 = \sum \frac{\omega_k}{2} (p_k^2+q_k^2).\] 
After the following canonical transformation to  complex  variables (analogs of  $\hat a$, $\hat a^\dag$  in quantum mechanics),
\begin{equation}\begin{cases} q_k &= \frac{1}{\sqrt{2}} (\xi_k + \rmi\, \eta_k),\\ p_k &= \frac{\rmi}{\sqrt{2}} (\xi_k - \rmi\, \eta_k),\end{cases} 
\label{zetaeta}\end{equation}
the Hamiltonian becomes
\[H=\sum_{k=1}^{\textbf d} \rmi\,\omega_k \xi_k \eta_k + \alpha \sum_{|\vecp m|+|\vecp n|\ge 3} \tilde H_i(\vecp m,\vecp n) \prod_{k=1}^{\textbf d} \xi_k^{m_k} \eta_k^{n_k}.
\]
The monomials $\xi^{\vecp m}\eta^{\vecp n}=\prod \xi_k^{m_k} \eta_k^{n_k}$  are  eigenvectors of the unperturbed Liouvillian 
\[\Lo_{H_0} \xi^{\vecp m}\eta^{\vecp n} = \rmi\,(\vecp\omega,\vecp m - \vecp n)\, \xi^{\vecp m}\eta^{\vecp n}.\]
 Therefore, for any series $F({\vecp p,\vecp q})= \sum \tilde F(\vecp m,\vecp n)  \xi^{\vecp m}\eta^{\vecp n}$ :
\begin{align}
\Po_{H_0} F &= \sum_{(\vecp\omega,\vecp m - \vecp n) =0} F(\vecp m,\vecp n)\, \xi^{\vecp m}\eta^{\vecp n},\nonumber\\
\So_{H_0} F &= \sum_{(\vecp\omega,\vecp m - \vecp n)\ne 0} \frac{1}{\rmi\,(\vecp\omega,\vecp m - \vecp n)}
F(\vecp m,\vecp n)\, \xi^{\vecp m}\eta^{\vecp n}.\nonumber
\end{align}
This is only a small excerpt from  the greater area of normal form theory \cite{Bruno,Sanders}.
\item{\bf ``Algebraic approaches.''} \cite{Lopatin,Bogaevski} These solve the homological equation using matrix methods in the enveloping algebra.
\item{\bf ``Zhuravlev quadrature.''} \cite{Zhuravlev02,petrov} This directly integrates $F(\bf x)$ along the unperturbed solution ${\bf x}(t,{\bf x_0})$. The asymptotic of the single quadrature
\[\int_0^{T}  F(\vecp p(t,\vecp q_0,\vecp p_0)),\vecp q(t,\vecp q_0,\vecp p_0)) \rmd t,\quad T\to\infty\]
 contains both $\Po_{H_0} F$ and $\So_{H_0} F$  inside its  $\Ord(T)$ and $\Ord(1)$ parts, respectively. 
This quadrature is very efficient in normal form theory. 
\end{itemize}

\subsection{Deprit perturbation series.}

If one chose $G_0=-\So_{H_0}H_i$ in $(\ref{Htrans})$, then all  terms of  order  $\alpha$  in the  transformed Hamiltonian will 
become secular (begin with $\Po_{H_0}$), as follows: 
\begin{equation*}
\widetilde H=\iexpD{\alpha\Lo_G} (H_0+\alpha H_i) 
= H_0+\alpha \Po H_i
+\frac{\alpha^2}{2}\left( \Lo_{H_0} G_1  + \Lo_{\So H_i}(1+\Po)H_i \right)
+\Ord(\alpha^3).
\end{equation*}
Hereafter, we will omit the subscript $H_0$ for   { unperturbed} $\Po_{H_0}$ and $\So_{H_0}$ operators. 

 In the next orders, one can consequentially choose $G_1= - \So \Lo_{\So H_i}(1+\Po)H_i$ 
to eliminate nonsecular terms up to $\alpha^2$, then  $G_2$ to eliminate nonsecular terms up to $\alpha^3$ and so on. 
We will refer to this process as  to the programme of  classical Poincar\'e-Lindstedt perturbation theory:
{\em using near-identity  canonical transformation, turn the Hamiltonian into an integral of the unperturbed system.}

Compare this to quantum mechanical perturbation theory:
{\em using near-identity  unitary transformation, turn the Hamiltonian operator into a block-diagonal form (commutative with the unperturbed Hamiltonian).}

The described procedure recursively constructs a generator of  normalizing transformation up to any order in $\alpha$. One may find details of  ``Deprit's triangular algorithm'' in  classical books on perturbation theory \cite{Naifeh,Giacaglia}. 
Here we will search for regularities in perturbation series.  Look at the first few orders:
\begin{multline}
G=-\So H_i - \alpha \So \Lo_{\So H_i}(1+\Po)H_i 
-\alpha^2\So  \left(\tfrac{1}{2} \Lo_{\So H_i} \Lo_{\Po H_i}\right.\\ 
\shoveright{{}+
( \left.\Lo_{H_i} \So - \Lo_{\So H_i} \Po + \tfrac{1}{2} \Lo_{\Po H_i} \So ) 
 \Lo_{(1+\Po) H_i} \right)\So H_i
+\Ord(\alpha^3)}\\
\shoveleft{\widetilde H
= H_0+\alpha \Po H_i
+\frac{\alpha^2}{2} \Po \Lo_{\So H_i}(1+\Po )H_i
+\frac{\alpha^3}{3} \Po  \left(\tfrac{1}{2} \Lo_{\So H_i} \Lo_{\Po H_i}\right.}
\label{Nonsimp}\\
{{}\left. +( \Lo_{H_i} \So - \Lo_{\So H_i} \Po + \tfrac{1}{2} \Lo_{\Po H_i} \So ) 
 \Lo_{(1+\Po) H_i} \right)\So H_i
+\Ord(\alpha^3)}
\end{multline}
At  first glance, there is no notable pattern in the above expressions. However, in 
the subsequent sections of this work, they will be transformed into a more systematic form using canonical identities.

Obviously, the  generator $G$ is not unique. An arbitrary secular function may be added to any order. Further orders  
will depend on this term. The rule for choosing the secular part of $G$ is called  {\em (hyper-)normalization  style} \cite{Sanders}.
Usually, one constructs a completely nonsecular  generator of Poincar\'e-Lindstedt transformation  
$\Po_{H_0} G = 0$.
However, we will see that other conditions may be useful. 

It is now necessary to say a few words  about the integrability of normalized system. In the non-resonant case, the unperturbed 
{\bf d}-dimensional system has only {\bf d}  integrals of motion. 
The normalized non-resonant Hamiltonian will be a function of these ``actions'' only,
and are formally integrable (up to $\Ord(\alpha^n))$.

However, if the unperturbed system has {\bf r} independent resonance relations  $(\vecp \omega,\vecp D_k)=0$, $k=1,\ldots ,{\textbf r}$,
then it has {\bf r} additional non-commutative integrals. 
There  still exists {\bf d}  involutive actions, but normalized Hamiltonian will also
 depend on phases.
 Even in this resonant case, one can decrease the order of perturbed system by at least  the dimension of center of  the algebra of unperturbed integrals. 
This center consists of ${\textbf d}-{\textbf r}$ commutative with all other integrals \cite{Gustavson}.
Reduction by these variables \cite{Arnold,Moser} results in an effective {\bf r} dimensional  system.
Therefore, in case that two or more
resonances relations exist, the normalized system will be generally non-integrable.

Here is the difference between classical and quantum-mechanical perturbation theories.
In quantum mechanics, one can always completely diagonalize a Hamiltonian  with any number of resonant
relations. Additional quantum integrals and the connection of  diagonalizing transformation to nonlinear 
finite-dimensional Bogolyubov transformations were discussed in \cite{Nikolaev,Nikolaev2}. 

\section{Canonical identities.}

In order to simplify expressions, we must list the algebraic properties of basic perturbation operators
$\Lo$, $\Po_{H_0}$ and $\So_{H_0}$:
\begin{itemize} 
\item  Invariant definitions  of $\Po$   and $\So$  operators $(\ref{Pdef},\ref{Sdef})$
immediately lead to  identities
\begin{align}
\Po H_0  \,\,&= H_0,\qquad \So H_0 = 0,  \nonumber\\
\Po\Lo_{H_0}&=\Lo_{H_0}\Po=0,\label{Fid}\\
\So\Lo_{H_0}&=\Lo_{H_0}\So = 1 - \Po.\nonumber
\end{align}
\item Because canonical transformation preserves algebraic and Poisson products (brackets), then 
for any functions $F({\bf x})$, $G({\bf x})$  and $H_i({\bf x})$ the following hold:
\begin{align}
\Po ( F \cdot\Po G )&=(\Po F) \cdot (\Po G),\nonumber\\
\Po\Lo_{\Po H_i}&=\Lo_{\Po H_i}\Po,\label{Pid}\\
\So ( F \cdot\Po G )&=(\So F)\cdot (\Po G),\nonumber\\
\So\Lo_{\Po H_i}&=\Lo_{\Po H_i}\So.\nonumber
\end{align}

The first two identities demonstrate that the projector
$\Po$  preserves the products if one operand already belongs to the algebra of integrals of motion.
 This is the canonical analogue of orthogonal projection.

\item Since $\Lo$ is the derivation for both algebraic and Poisson products of functions, there  exists the {\em integration by parts} formulae
known as { Friedrichs identity} \cite{Friedrichs}:
\begin{align}
\So\Lo_{H_i}\So &=\Lo_{\Po H_i}\So^2+ \Lo_{\So H_i}\So - \So\Lo_{\So H_i}  -\Po\,\Lo_{\So H_i}\So  + \So\,\Lo_{\So H_i}\Po,\label{Friedrichs}\\
\So (F\cdot\So G)&=\Po F\cdot \So^2 G + ( \So F\cdot\So G) - \So ( \So F\cdot G) - \Po ( \So F\cdot\So G) + \So ( \So F\cdot\,\Po G)\nonumber
\end{align}
The proof follows from an application
of both sides of
$\So \Lo_{H_0} = 1-\Po$ to the product $\Lo_{\So H_i}\So F=[\So F,\So H_i]$ (or $\So F\cdot \So G$),  and an  
expansion of the Jacobi identity.  

 We may achieve  complete  symmetry  denoting the algebraic product as an operator. Hereafter, we will  even omit the identities for algebraic
product, automatically assuming the existence  of complementary identities. 

\item Further identities exist for products of three basic operators. For any function $H_i$:
\begin{align}
\Po\Lo_{H_i}\Po&=\Lo_{\Po H_i}\Po,\label{Eid}\\
\So\Lo_{H_i}\Po&=\Lo_{\So H_i}\Po.\nonumber
\end{align}
These  are the consequences of Poisson bracket invariance under canonical transformations
\begin{equation*}
\Po\Lo_{H_i}\Po F = \lim_{T\to\infty}\frac{1}{T} \int_0^{T} \rmd t\, \rme^{t\Lo_{H_0}}[ \Po F, H_i ] 
= \lim_{T\to\infty}\frac{1}{T} \int_0^{T} \rmd t\,  [ \rme^{t\Lo_{H_0}} \Po F, \rme^{t\Lo_{H_0}} H_i ]. 
\end{equation*}
Since $\Po F$ is the integral of motion for $H_0$,
\[
\Po\Lo_{H_i}\Po F= [\Po F, \lim_{T\to\infty}\frac{1}{T} \int_0^{T} \rmd t\, \rme^{t\Lo_{H_0}} H_i]=\Lo_{\Po H_i}\Po F.
\]
The proof of the second identity is similar.
\item The Burshtein-Soloviev identity \cite{Burshtein}
\begin{equation}
\Po\Lo_{H_i}\So =-\Po\Lo_{\So H_i},\label{Lid}
\end{equation}
 follows from identically zero expression $\Po\Lo_{H_0}\Lo_{\So H_i}\So\equiv 0$. Using the Jacobi identity
$\Lo_F \Lo_G - \Lo_G \Lo_F = \Lo_{\Lo_F G}$, we can write
\small\[
0  \equiv
\Po\Lo_{H_0}\Lo_{\So H_i}\So =\Po\Lo_{\So H_i}\Lo_{H_0}\So + \Po\Lo_{\Lo_{H_0}\So H_i}\So 
=\Po\Lo_{\So H_i} + \Po\Lo_{H_i}\So - \Po\Lo_{\So H_i}\Po - \Po\Lo_{\Po H_i}\So.
\]\normalsize
Here the last two terms  vanish, due to (\ref{Eid}) and (\ref{Pid}). Therefore,
$ \Po\Lo_{\So H_i} + \Po\Lo_{H_i}\So \equiv 0$.
\end{itemize}
The identities $(\ref{Fid})-(\ref{Lid})$ is all that we need to simplify the Deprit series.  They could be directly verified in an ``Action-Angle'' representation.
Later we will find their generalization. 

We have used the computer algebra system FORM \cite{Form} to implement the above formulae.
Due to  canonical identities, the first orders of the Deprit series $(\ref{Nonsimp})$  have been reduced to very inspiring form
(see also the   ``{\tt deprit\_series}'' demo):
\begin{multline}
G=-\So H_i+\alpha\left(\So\Lo\So H_i- \So^2\Lo\Po H_i\right)
+\alpha^2 \left(-\So\Lo\So\Lo\So H_i   \right.\\\left.
\shoveright{+  \So\Lo\So^2\Lo\Po H_i+\So^2\Lo\So\Lo\Po H_i
             + \So^2\Lo\Po\Lo\So H_i\right)+ \Ord(\alpha^3)}\\
\shoveleft{\tilde H=H_0 
       + \alpha \Po H_i 
       - \frac{\alpha^2}{2} \Po\Lo\So H_i 
       + \alpha^3 \left(\tfrac{1}{3}\Po\Lo\So\Lo\So H_i
                - \tfrac{1}{6}\Po\Lo\So^2\Lo\Po H_i \right)+\Ord(\alpha^4)}\label{Simp}
\end{multline}

The expression for the transformed Hamiltonian corresponds to the classic result of Burshtein and Soloviev \cite{Burshtein,mitropolsky71}.   In  standard notation, 
this is
\[\widetilde H=H_0 + \alpha\overline H_i+\frac{\alpha^2}{2}\overline{\left[ \widetilde{H_i}, H_i \right]} + \frac{\alpha^3}{3}\overline{\left[ \widetilde{H_i}, \left[\widetilde{H_i}, H_i + \frac{\overline{H_i}}{2}\right] \right]}+\Ord(\alpha^4).\]

It is intriguing that the expression for $G$ looks like a simple sum of all compositions of the $-\So$ and $\Po$ operators. 
 In the third order, the  expression for the non-secular style Deprit generator looses this structure.
Nevertheless, we  observed that such sums actually normalize the Hamiltonian in the next orders.
This should be explored further.
\section{Kato expansion.}
\subsection{Resolvent of Liouville operator}
Quantum mechanics commonly uses the stronger  {\em Abel  averaging}  \cite{HP} procedure
\[
\langle F\rangle^{(A)} =
\lim_{\lambda\to +0}\lambda\int_0^{+\infty}\rme^{-\lambda t} \rme^{t\Lo_{H_0}} F({\bf x})\, \rmd t.
\]
This can  also be applied to the classical case. Corresponding averaging and integrating operators are known from quantum mechanics  \cite{Primas,Tyuterev,Jauslin}:
\begin{align}
&\Po_{H_0} = \lim_{\lambda\to +0}\lambda\int_0^{+\infty}\rmd t \,\rme^{-\lambda t} \rme^{t\,\Lo_{H_0}},\label{Adef}\\
&\So_{H_0} = - \lim_{\lambda\to{+0}}\int_0^\infty \rmd t \,\rme^{-\lambda t} \rme^{t\,\Lo_{H_0}} \left(1-\Po_{H_0}\right).\nonumber
\end{align}
Whenever the { Ces\`aro}  average $(\ref{Cesaro})$ exists,   { Abel} averaging gives the same results \cite{HP}. This is why we use the same notation. 
From this point, we will always assume  Abel averaging, if not  otherwise specified, 
as this greatly simplifies formulae  and forms a natural  connection to resolvent formalism. Strictly speaking, we should discuss corresponding Tauberian theorems, 
but we have limited our goal   to  formal expressions only.  

The Abel averaging definitions for $\So_{H_0}$  and $\Po_{H_0}$,  as well as quantum mechanical analogies, 
 suggest exploring the resolvent of Liouville operator 
\begin{equation}
\R_{H}(z) =\frac{1}{{\Lo_{H}-z}}.\label{Res}
\end{equation}
This  operator-valued function of the complex variable $z$ is the Laplace transform of  the evolution operator   of Hamiltonian system
\[
\R_{H}(z) = - \int_0^{+\infty}\rmd t\, \rme^{-z t} \rme^{t\Lo_{H}}.
\]

Resolvent singularities are the eigenvalues of   $\Lo_{H}$. For integrable Hamiltonian system with compact energy surfaces, these eigenvalues belong to an imaginary axis. 
Typically, the  spectrum of a Liouville operator is anywhere dense \cite{spohn75}. 

But let us begin with the simpler case of an isolated point spectrum. We will consider one-dimensional  system and 
 restrict the  domain of resolvent operator to  analytic functions 
with argument {\it on the compact energy surface} $H({\bf x})=E$.
Under such conditions the system is non-relaxing and oscillates with the single frequency $\omega(E)$. The resolvent singularities are located  at points 0, $\pm \rmi\omega(E)$, $\pm 2\rmi\omega(E)$, \ldots

We are interested in the analytical structure of the resolvent around zero. The existence of  $\Po_{H_0}$ and $\So_{H_0}$ $(\ref{Adef})$
means that the unperturbed resolvent has a simple pole in $0$. 
The averaging operator is the residue of the resolvent in this pole
\[
\Po_{H_0} \equiv - \res_{z=0} \R_{H_0},
\]
while the integrating operator $\So_{H_0}$ is its holomorphic part 
\[
\So_{H_0}=\lim_{z\to0}{\R_{H_0}(z)(1-\Po_{H_0})}.
\] 
Therefore, the Liouvillian resolvent combines both basic perturbation operators  \cite{nikolaev3}.

We will also need the  Hilbert identity, 
\begin{equation}
\R_{H}(z_1)-\R_{H}(z_2)=(z_1-z_2)\,\R_{H}(z_1)\R_{H}(z_2),\label{HilbertId}
\end{equation}
which holds for any  complex $z_1$ and $z_2$ outside of the spectrum (resolvent set) of $\Lo_H$.

\paragraph{The Laurent series.}
Despite the fact that the unperturbed resolvent has the simple pole in the origin, the perturbed resolvent  may be more singular. The Laurent series of
the  resolvent with an isolated singularity in origin is
\begin{equation}
\R_{H}(z)=\sum_{n=-\infty}^{+\infty} \Sr_{H}^{(n)}z^{n-1},\quad \Sr_H^{(n)} = \frac{1}{2\pi \rmi}\oint_{|z|=\epsilon}\R_H(z)\,z^{-n}\,\rmd z \label{u1}
\end{equation}
(here $(n)$ is index).

It is well known that the Hilbert resolvent identity $(\ref{HilbertId})$  defines the structure of this series \cite{Kato}. 
For further usage, it makes sense to repeat here this derivation.
Consider the product of the two coefficients 
\begin{multline}
\Sr_H^{(m)}\,\Sr_H^{(n)} = \left(\frac{1}{2\pi \rmi}\right)^2
\oint_{|z_1|=\epsilon_1}\oint_{|z_2|=\epsilon_2}
\R_H(z_1)\R_H(z_2)z_1^{-m}\,z_2^{-n}\,\rmd z_1\,\rmd z_2\\
= \left(\frac{1}{2\pi \rmi}\right)^2
\oint_{|z_1|=\epsilon_1}\oint_{|z_2|=\epsilon_2}
\frac{z_1^{-m}\,z_2^{-n}}{ z_2-z_1}\left(\R_H(z_2)-\R_H(z_1)\right)\label{ResCoeff}
\,\rmd z_1\,\rmd z_2.
\end{multline}
Integration around the two circles $\epsilon_1 < \epsilon_2 < \omega$ gives
\begin{align*}
{1\over2\pi \rmi} \oint_{|z_1|=\epsilon_1} \frac{z_1^{-m}}{ z_2 - z_1} \rmd z_1
&=\eta_m z_2^{-m},\qquad{\rm where}\,\,\eta_n = \begin{cases}1&{\rm when}\, n\ge 1,\\0&{\rm when}\, n<1,\end{cases}\\
{1\over2\pi \rmi} \oint_{|z_2|=\epsilon_2} \frac{z_2^{-n}}{ z_2 - z_1} \rmd z_2
&=(1-\eta_n) z_1^{-n}.&
\end{align*}

As a result, the resolvent operator coefficients obey 
\begin{displaymath}
\Sr_H^{(m)}\,\Sr_H^{(n)} = (\eta_m + \eta_n -1) \Sr_H^{(m+n)}.
\end{displaymath}
For $n=m=0$, this states that resolvent residue (with a minus sign) is the projector
\[
\Sr_{H}^{(0)}= - \Sr_{H}^{(0)} \Sr_{H}^{(0)} = -\Po_{H},
\]
and 
\begin{align*}
&\Sr_H^{(n)}=\So^n_H,\qquad n\ge 1,\\
&\Sr_H^{(-n)}\!\!=-\D^n_H,\\
&\So_{H}\Po_{H}=\Po_{H}\So_{H}\equiv 0,\quad \So_{H}\D_{H}=\D_{H}\So_{H}\equiv 0,\\ 
&\Po_{H}\D_{H}=\D_{H}\Po_{H} = \D_{H}.
\end{align*}
Here, $\D_H$ is the eigennilpotent operator, which does not have an unperturbed analogue ($\D_{H_0} \equiv 0$). 
Therefore,  the Laurent series of the general Liouvillian resolvent  around $z=0$ has the form \cite{Kato}
\begin{equation}
\R_H(z)=-{1\over z}\Po_H + \sum_{n=0}^\infty z^n \So_H^{n+1}
-\sum_{n=2}^\infty z^{-n} \D_H^{n-1},\label{ResLoran}
\end{equation}
while the unperturbed resolvent consists  only of
\[
\R_{H_0(z)}=\sum_{n=0}^{+\infty} \Sr_{H_0}^{(n)}z^{n-1}=-{1\over z}\Po + \sum_{n=0}^\infty z^n \So^{n+1}\label{UResLoran}
\]
(remember that we omit the subscript $H_0$).

\paragraph{Canonical structure of resolvent.}
Yet another  identity of the Hilbert type relates the  resolvent of the Liouville operator to the Poisson brackets.
For any $z_1$, $z_2$, $z_3$  outside of the spectrum  of $\Lo_H$ and functions $F$ and $G$ the following hold true:
\begin{multline}
\R_{H}(z_1)[\R_{H}(z_2)F,G] + \R_{H}(z_1)[F,\R_{H}(z_3)G]
-[\R_{H}(z_2)F,\R_{H}(z_3)G]=\\
=(z_1-z_2-z_3)\R_{H}(z_1)[\R_{H}(z_2)F,\R_{H}(z_3)G].\label{SimpID}
\end{multline}
This is another  {\em integration by parts} formula. It follows from the application  of  identical operator
$
\R_{H}(z_1)(\Lo_H - z_1) \equiv  1
$
to Poisson bracket $[\R_{H}(z_2)F,\R_{H}(z_3)G]$ 
and the expansion of the Jacobi identity. 
A complementary identity exists for the algebraic product. 

We will use  this canonical identity in the operator form
\begin{multline}
\R_H(z_1) \Lo_{\R_H(z_2) F} - \Lo_{\R_H(z_2) F} \R_H(z_3)  + \R_H(z_1) \Lo_{F}  \R_H(z_3) =\\
= (z_1-z_2-z_3) \R_H(z_1) \Lo_{\R_H(z_2) F}  \R_H(z_3).\label{SimpID1}
\end{multline}
 
Its Laurent coefficients for the unperturbed resolvent contain all previous canonical identities $(\ref{Fid})-(\ref{Lid})$.  
 For example, the coefficient of $z_1^0 z_2^0 z_3^0$ is the Friedrichs identity $(\ref{Friedrichs})$, and so on.  
Integrating $(\ref{SimpID1})$ like
  Hilbert identity before, we can obtain advanced  identities. 
We will do this for a perturbed resolvent in the next chapter.

\subsection{Kato series}
The perturbed resolvent can be expanded into  the Neumann series as follows:
\begin{multline*}
\R_{H_0+\alpha H_i}(z)=
\R_{H_0}-\alpha \R_{H_0}\Lo_{H_i}\R_{H_0}+
\alpha^2 \R_{H_0}\Lo_{H_i}\R_{H_0}\Lo_{H_i}\R_{H_0}+
\ldots\\
=\sum_{n=0}^{\infty}{(-1)}^n\alpha^n\R_{H_0}(z){\left(\Lo_{H_i}
\R_{H_0}(z)\right)}^n.
\end{multline*}


The integration 
around a small contour results in the Kato   series \cite{Kato} for the ``perturbed averaging operator''
\begin{multline*}
\Po_H = -{1\over{2\pi \rmi}}\oint_{|z|=\epsilon} \R_H(z)\rmd z
=-{1\over{2\pi \rmi}}\sum_{n=0}^\infty \oint_{|z|=\epsilon} {(-1)}^n\alpha^n
\R_{H_0}(z){\left(\Lo_{H_i}\R_{H_0}(z)\right)}^n \, \rmd z\\
={-1\over{2\pi \rmi}}\sum^\infty_{n=0} {(-1)}^n\alpha^n\!\oint_{|z|=\epsilon}
\left(\sum_{m=0}^\infty \Sr_{H_0}^{(m)}z^{m-1}\right)
{\left(\Lo_{H_i}\sum_{k=0}^\infty\Sr_{H_0}^{(k)}z^{k-1}\right)}^n \rmd z ,
\end{multline*}
 the ``perturbed integrating operator'', and the ``perturbed quasi-nilpotent'':
\begin{equation*}
\So_H = {1\over{2\pi \rmi}}\oint_{|z|=\epsilon} z^{-1}\R_H(z)\rmd z,\qquad \D_H = -{1\over{2\pi \rmi}}\oint_{|z|=\epsilon} z\R_H(z)\rmd z .
\end{equation*}
Only  coefficients of $z^{-1}$ in these expansions will contribute to the result:
\begin{align}
\Po_H&=\sum_{n=0}^\infty{(-1)}^{n+1}\alpha^n\left(\sum_{
\sum{p_i}={\bf n}\atop p_i\ge 0}\Sr_{H_0}^{(p_{n+1})}\underbrace{\Lo_{H_i}\Sr_{H_0}^{(p_n)}\,\ldots\,\Sr_{H_0}^{(p_2)}
\Lo_{H_i}}_{n\ {\rm times}}\Sr_{H_0}^{(p_1)}\right),\nonumber\\
\So_H&=\sum_{n=0}^\infty{(-1)}^{n\hphantom{+1}}\alpha^n\left(\sum_{
\sum{p_i}={\bf n+1}\atop p_i\ge 0}\Sr_{H_0}^{(p_{n+1})}\underbrace{\Lo_{H_i}\Sr_{H_0}^{(p_n)}\,\ldots\,\Sr_{H_0}^{(p_2)}
\Lo_{H_i}}_{n\ {\rm times}}\Sr_{H_0}^{(p_1)}\right),\label{Pexp}\\
\D_H&=\sum_{n=1}^\infty{(-1)}^{n+1}\alpha^n\left(\sum_{
\sum{p_i}={\bf n-1}\atop p_i\ge 0}\Sr_{H_0}^{(p_{n+1})}\underbrace{\Lo_{H_i}\Sr_{H_0}^{(p_n)}\,\ldots\,\Sr_{H_0}^{(p_2)}
\Lo_{H_i}}_{n\ {\rm times}}\Sr_{H_0}^{(p_1)}\right).\nonumber
\end{align}
A summation in the above expressions should be done by all possible placements  of $n$ (or $n+1$, or  $n-1$)  operators $\So_{H_0}$   in $n+1$ sets. 
These are also known as a ``weak compositions''. There are $C^n_{2n}$ terms of order $\alpha^n$  in $\Po_H$  and $C^n_{2n+1}$ terms in $\So_H$.
Here are the first two orders of the ``perturbed integrating operator'':
\begin{multline*}
\So_H=\So -\alpha(\So\Lo_{H_i}\So - \So^2\Lo_{H_i}\Po - \Po\Lo_{H_i}\So^2) 
+\alpha^2(\So\Lo_{H_i}\So\Lo_{H_i}\So-\So^2\Lo_{H_i}\So\Lo_{H_i}\Po\\
{}-\! \So^2\Lo_{H_i}\Po\Lo_{H_i}\So
\! -\! \Po\Lo_{H_i}\So^2\Lo_{H_i}\So\! -\! \So\Lo_{H_i}\So^2\Lo_{H_i}\Po\!  -\! \So\Lo_{H_i}\Po\Lo_{H_i}\So^2\!\!-\Po\Lo_{H_i}\So\Lo_{H_i}\So^2\\
{}+\Po\Lo_{H_i}\Po\Lo_{H_i}\So^3
+\Po\Lo_{H_i}\So^3\Lo_{H_i}\Po +\So^3\Lo_{H_i}\Po\Lo_{H_i}\Po) + \Ord(\alpha^3),
\end{multline*}
the ``perturbed  projector'':
\begin{multline*}
\Po_H=\Po -\alpha(\Po\Lo_{H_i}\So + \So\Lo_{H_i}\Po)
+\alpha^2(\Po\Lo_{H_i}\So\Lo_{H_i}\So +\So\Lo_{H_i}\Po\Lo_{H_i}\So \\
{}+ \So\Lo_{H_i}\So\Lo_{H_i}\Po -\Po\Lo_{H_i}\Po\Lo_{H_i}\So^2 
- \Po\Lo_{H_i}\So^2\Lo_{H_i}\Po 
- \So^2\Lo_{H_i}\Po\Lo_{H_i}\Po)+ \Ord(\alpha^3)
\end{multline*}
and, finally, the ``perturbed quasi-nilpotent'':
\begin{multline*}
\D_H=\alpha\Po\Lo_{H_i}\Po -\alpha^2(\Po\Lo_{H_i}\Po\Lo_{H_i}\So +
\Po\Lo_{H_i}\So\Lo_{H_i}\Po +\So\Lo_{H_i}\Po\Lo_{H_i}\Po)
+\Ord(\alpha^3).
\end{multline*}\Checked

Properties of unperturbed operators can be  extended to their analytic continuations 
as follows:
\[\Po_H H = H, \quad\So_H\Lo_H = 1-\Po_H,\quad\Lo_H \Po_H = \Po_H \Lo_H = \D_H, \, {\textrm{etc.}}\]
For  details, see  Appendix  A and the demo   ``{\tt perturbed\_operators}''. 

Quantum mechanical perturbation theory uses series for $\Po_H$  in order to find  the eigenvalues of perturbed Hamiltonian \cite{messiah}. 
There is no straightforward analogue  in classical mechanics.
To avoid misunderstanding, it should be noted that $\Po_H F$ will not be an integral  of the perturbed Hamiltonian. 
This is because $\Lo_H \Po_H =  \D_H$  is nonzero, in general. 

Actually, the  ``perturbed projector''  $\Po_H$ projects onto the analytic continuation
of the algebra of integrals of unperturbed Hamiltonian. 
This does not coincide, in general, with algebra of integrals of perturbed system.
 In other words, the zero eigenvalue of the Liouville operator may be split by perturbation.

\subsection{Connection to Poincar\'e-Lindstedt-Deprit series}

To establish a relation between  Kato expansion  and standard classical perturbation theory,  
we will demonstrate that the projectors $\Po_H$ and $\Po_{H_0}$ are connected by canonical transformation.
This is the junction point of  two very different formalisms. 

Indeed, if we   construct  a transformation
\begin{displaymath}
\tilde {\bf{x}} = \expD{\alpha\Lo_G} \, {\bf{x}},\qquad \widetilde H = \iexpD{\alpha\Lo_G}H,
\end{displaymath}
connecting the projectors
\begin{equation}
\Po_H = \expD{\alpha\Lo_G} \Po_{H_0} \iexpD{\alpha\Lo_G},\label{Symeq}
\end{equation}
then it follows  from  $\Po_H H = H$ that the transformed Hamiltonian will be  an integral of the unperturbed system 
\begin{equation*}
\Po_{H_0} \widetilde H = \Po_{H_0} \iexpD{\alpha\Lo_G} H = 
\iexpD{\alpha\Lo_G} \Po_H H 
= \iexpD{\alpha\Lo_G} H = \widetilde H .
\end{equation*}
Therefore, such transformation will directly realize the programme of Poincar\'e-Lind\-stedt perturbation theory.

To find the transformation, consider the derivative of the resolvent  with respect to the perturbation
\[
\frac{\partial}{\partial\alpha}\R_H(z) = -\R_H(z)\Lo_{H_i}\R_H(z).
\]
Substitution $z_1=z_3=z$ and $F=H_i$  into the canonical resolvent identity $(\ref{SimpID1})$  results in
\small\begin{equation}
\R_H(z)\Lo_{H_i}\R_H(z) = \Lo_{\R_H(z_2) H_i}\R_H(z) - \R_H(z)\Lo_{\R_H(z_2) H_i} - z_2  \R_H(z) \Lo_{\R_H(z_2) H_i}\R_H(z) \label{SimpID2}
\end{equation}\normalsize
Look at the coefficient of $z_2^0$ in the Laurent series of the above  expression 
\[
\frac{\partial}{\partial\alpha}\R_H(z) = 
\R_H(z)\Lo_{\So_H H_i} - \Lo_{\So_H H_i}\R_H(z) - \R_H(z)\Lo_{\Po_H H_i}\R_H(z).
\]
Proceeding the same way for coefficients of $z_2^{-n}$ \hbox{$(n\ge1)$} in $(\ref{SimpID2})$, 
we find 
\begin{align}
\R_H(z)\Lo_{\Po_H H_i}&=\Lo_{\Po_H H_i}\R_H(z) - \R_H(z)\Lo_{\D_H H_i}\R_H(z),
\nonumber\\
\R_H(z)\Lo_{\D^n_H H_i}&=\Lo_{\D^n_H H_i}\R_H(z) - \R_H(z)\Lo_{\D^{n+1}_H H_i}\R_H(z).
\nonumber
\end{align}\Checked
This allows for the rewriting of the  resolvent derivative as
\begin{align}
{\partial\over \partial\alpha}\R_H(z) &=
\R_H(z)\Lo_{\So_H H_i} - \Lo_{\So_H H_i}\R_H(z) 
\nonumber\\
&- \Lo_{\Po_H H_i}\R_H(z)^2 
+ \Lo_{\D_H H_i}\R_H(z)^3
- \Lo_{\D^2_H H_i}\R_H(z)^4 + \ldots\nonumber
\end{align}
\Checked
Actually, this is a power series because $\D^n_H = \Ord(\alpha^n)$.

From the Hilbert identity, it follows that
${\partial^n\over \partial z^n}\R_H(z) = n!\R_H^{n+1}(z)$.
Finally,
\begin{multline}
\frac{\partial}{\partial\alpha}\R_H(z) =
\R_H(z)\Lo_{\So_H H_i} - \Lo_{\So_H H_i}\R_H(z) 
- \Lo_{\Po_H H_i}\frac{\partial\,\R_H(z)}{\partial\,z} \\
{}+ \frac{1}{2} \Lo_{\D_H H_i}\frac{\partial^2\R_H(z)}{\partial\,z^2}
- \frac{1}{6} \Lo_{\D_H^2 H_i}\frac{\partial^3\R_H(z)}{\partial\,z^3} 
+ \ldots
\end{multline}
\Checked
Therefore, change of the resolvent under perturbation $\alpha H_i$ can be represented as a sum of canonical transformation with generator  $-\So_H H_i$
and resolvent transformation as a function of complex variable $z$.

The derivative of projector $\frac{\partial}{\partial\,\alpha}\Po_H$ may be obtained as the residue of the previous expression at $z=0$.
In our  case of isolated point spectrum,  the resolvent  is a meromorphic function, and the residue of any its derivative with respect to $z$ vanishes. 
Therefore, the projector $\Po_H$  transforms canonically under  perturbation
\begin{equation}
\frac{\partial}{\partial\,\alpha}\Po_H = \Po_H\Lo_{\So_H H_i} - \Lo_{\So_H H_i}\Po_H,\label{Projeq}
\end{equation}
and the projectors are connected by the Lie transform 
with generator $-\So_H H_i$:
\[
\Po_H = {\expD{\alpha\Lo_{-\So_H H_i}} \Po_{H_0} \iexpD{\alpha\Lo_{-\So_H H_i}}}.
\]
We see   that  the canonical transformation with generator
\begin{multline}
G=-\So_H H_i=\sum_{n=0}^\infty{(-1)}^{n+1}\alpha^n\left(\sum_{
\sum
{p_i}={\bf n+1}\atop p_i\ge 0}\Sr_{H_0}^{(p_{n+1})}
\underbrace{\Lo_{H_i}\Sr_{H_0}^{(p_n)}\,\ldots\,\Sr_{H_0}^{(p_2)}
\Lo_{H_i}}_{n\ {\rm times}}\Sr_{H_0}^{(p_1)} H_i
\right) \label{explicit}
\end{multline}
formally normalizes the Hamiltonian in all orders in $\alpha$.

\subsection{General form of generator.}
Knowing that $\Po_{H_0}$ and $\Po_H$ are canonically connected  and that this transformation normalizes the Hamiltonian, 
we can   reformulate the programme of Poincar\'e-Lind\-stedt perturbation theory to 
{\em the construction of  canonical transformation, which connects unperturbed and perturbed averaging operators.}

Let us now determine the general form of
such  a transformation. It follows from  $(\ref{Symeq})$ that $\Po_H$ satisfies the operatorial differential equation
\[
\frac{\partial}{\partial\,\alpha}\Po_H = \Lo_G \Po_H - \Po_H \Lo_G.
\]
Application of  this expression to Hamiltonian $H$ results in
\begin{displaymath}
({\partial\over \partial\,\alpha}\Po_H) H  = \Lo_G \Po_H H - \Po_H \Lo_G H.\nonumber
\end{displaymath}
Since $\Po_H H = H$,  
${\partial\over \partial\,\alpha}(\Po_H H) = \left({\partial\over \partial\,\alpha}\Po_H\right) H + 
\Po_H \frac{\partial}{\partial\,\alpha}H $ and ${\partial\over \partial\,\alpha}H = H_i$,
the previous expression becomes
\[
(1-\Po_H) \Lo_H G = - (1-\Po_H) H_i.
\]
To solve this equation, it is sufficient to apply the $\So_H$ operator. Therefore, the general form of the generator of connecting transformation is
\begin{equation}
G = - \So_H H_i + \Po_H F,\label{Final}
\end{equation}\Checked
where $F({\bf x})$  may be any analytic function. This is the central formula of this work. It provides a non-recursive expression for
the generator of Poincar\'e-Lindstedt transformation and defines its ambiguity. 

The  choice of  function $F$ is the normalization {\em style}. It is natural to choose $F({\bf x})\equiv 0$ or $\Po_H G = 0$.
This is not equal to the ``non-secular'' normalization $\Po_{H_0} G_D = 0$  traditionally used for Deprit generator in classical perturbation theory. 
Because $\Lo_{\Po_H F} \Po_H = \Po_H \Lo_{\Po_H F}$, the projector $\Po_H$ itself is not sensitive to normalization style.

We may conclude that {\em the generators of normalizing transformations 
may differ by a function belonging to  a continuation of  
algebra of integrals of unperturbed system.}
Illustration can be found in  the demo   ``{\tt styles}'' in the Supplementary data files.

Quantum mechanical  Kato  perturbation expansion uses
another  formulae for 
unitary transformation connecting the unperturbed and perturbed projectors \cite{Kato}.
However, the original Kato generator   
and the rational expression developed by Sz\"okefalvi-Nagy \cite{Kato} 
do not define  canonical transformations  in classical mechanics.

\subsection{The first orders.}
Let us compare the  expressions for  generator
\begin{multline}
G=-\So H_i+\alpha\left(\So\Lo\So H_i- \So^2\Lo\Po H_i\right)
     - \alpha^2 (\So\Lo\So\Lo\So H_i- \So\Lo\So^2\Lo\Po H_i\\
{} -\So^2\Lo\So\Lo\Po H_i
             - \So^2\Lo\Po\Lo\So H_i-\Po\Lo\So\Lo\So^2 H_i- \Po\Lo\So^2\Lo\So H_i\\
{}  +\So^3\Lo\Po\Lo\Po H_i 
                 + \Po\Lo\So^3\Lo\Po H_i + \Po\Lo\Po\Lo\So^3 H_i) + \Ord(\alpha^3),
\end{multline}
 and  normalized Hamiltonian
\begin{multline}
\tilde H=H_0 
       + \alpha \Po H_i 
       -  \frac{\alpha^2}{2} \Po\Lo\So H_i 
       + \alpha^3 (\frac{1}{3} \Po\Lo\So\Lo\So H_i
                 - \frac{1}{6} \Po\Lo\So^2\Lo\Po H_i )
      {}+\alpha^4 (\frac{1}{6}\Po\Lo\So\Lo\So^2\Lo\Po H_i\\
                 {}- \frac{1}{4}  \Po\Lo\So\Lo\So\Lo\So   H_i 
                 + \frac{1}{12} \Po\Lo\So^2\Lo\So\Lo\Po H_i
               + \frac{1}{8} \Po\Lo\So^2\Lo\Po\Lo\So H_i  
               {}+ \frac{1}{4} \Po\Lo\Po\Lo\So^2\Lo\So H_i \\
                 + \frac{1}{4} \Po\Lo\Po\Lo\So\Lo\So^2 H_i 
                 - \frac{1}{6} \Po\Lo\Po\Lo\So^3\Lo\Po H_i 
  {}               - \frac{1}{4} \Po\Lo\Po\Lo\Po\Lo\So^3 H_i )+ \Ord(\alpha^5)\label{KatoGen}
\end{multline}
with  the Deprit series $(\ref{Simp})$.
Here we have denoted $\Lo=\Lo_{H_i}$  for compactness and used the identities $\Po\Lo_{H_i}\So^2 H_i \equiv 0$ and $\Po\Lo\Po H_i\equiv 0$.

At these  orders,  series \eref{Simp} and \eref{KatoGen} are very close and  differ  
only by secular terms in the generator.
This  is because of the natural normalization style $\Po_H G = 0$.
In  higher orders, the differences also propagate in non-secular terms. 
 Due to these additional terms, our expressions for the generator and normalized Hamiltonian become linear and systematic.
Actually, because of  canonical identities, 
there are many equivalent  expressions for generator and Hamiltonian. 
In  next section we will develop another explicitly secular one. 

Larger formula for  the Hamiltonian normalized up to $\Ord(\alpha^8)$  can be found in the demo ``{\tt normalized\_hamiltonian7}'' in the supplementary data files.
It consist of 528 terms and the corresponding generator has 2353 terms. 

\paragraph{Example 1.  Duffing equation.} This is  
 one-dimensional oscillator with the quartic anharmonicity:
\[H=\frac{1}{2}(p^2+q^2)+\frac{\alpha}{4} q^4=J+ \alpha J^2 \cos^4(\phi).
\]
The following are  Birkhoff and  ``Action-Angle'' representations of normalizing generator: 
\begin{align}
G&=\frac{p q}{32}   (3 p^2  + 5 q^2) -\alpha\frac{p q}{384}   (39 p^4 + 104 p^2 q^2 + 57 q^4)+\Ord(\alpha^2)\label{q4}\\
&=-{\frac{J^2}{32}\left( 8\sin (2\phi )+ \sin (4\phi ) \right) } +  {\frac{\alpha J^3}{192} \left( 99\sin (2\phi ) + 9\sin (4\phi )\! -\!  \sin (6\phi ) \right) }\!+\!\Ord(\alpha^2)\nonumber
\end{align}
The normalized Hamiltonian is as follows:
\begin{align*}
\tilde H&= \frac{1}{2} (p^2+q^2) + \frac{3\alpha}{32} {(p^2+q^2)}^2 -  \frac{17\alpha^2}{512} {(p^2+q^2)}^3 +  \frac{375 \alpha^3}{16384}{(p^2+q^2)}^4
+O(\alpha^4)\\
&=J + \frac{3}{8}\alpha J^2 -  \frac{17}{64}\alpha^2 J^3 +   \frac{375}{1024}\alpha^3 J^4 + O(\alpha^4).
\end{align*}
This is the typical structure of the perturbation series for nonresonant systems. More terms can be found in the demo  ``{\tt anharmonic}''.
Because of the uniqueness of the nonresonant normal form  \cite{Koseleff}, the normalized Hamiltonian coincides with  the classical Deprit series. 

\section{Multidimensional systems.}
\subsection{Truncated series.}
An extension of the previous construction to the multidimensional case causes  difficulty. 
The spectrum of the Liouville operator  
 for multi-frequency dynamical system 
is a union of countably many additive groups \cite{spohn75}. Typically, singularities of the multidimensional resolvent    cannot be separated from the origin. 
This is the classic problem of {\em small denominators}.
The resolvent  is no more holomorphic   \cite{Balescu75} and we can therefore no longer directly rely on Kato series \eref{Pexp} and  residues.

But our computer algebra calculations  confirm  the validity of the expression \eref{explicit} for generator 
  in the first orders of multidimensional case also.
The key point here is the canonical connection between  perturbed and unperturbed  projectors. 
Using  the demo ``{\tt phdot}'', we  checked  it by brute force  computation   up to 
$\alpha^{10}$. The simplification of almost three  million terms took a week of computer time.

Therefore, we can suggest that  formula \eref{explicit} remains (asymptotically) valid. 
But more proof is needed that  uses the canonical identities $(\ref{Fid})-(\ref{Eid})$ and does not rely on the resolvent and spectrum of Liouville operator. 

In the Appendices, we present this alternate proof. It   follows the ideas of the previous sections, 
but  uses  truncated at order $N$ finite  operatorial sums. 
Appendix A discusses the properties of  truncated operators $\Po_H^{[N]}$, $\So_H^{[N]}$ and $\D_H^{[N]}$, 
and in Appendix B we demonstrate that for any $N\in\mathbb{N}$, $N\ge2$, the  perturbed truncated projector transforms canonically:
\begin{equation}
 \left( \frac{\partial}{\partial\alpha}{\Po_H}\right)^{[N]} =    \Lo_{-\So_H^{[N]}  H_i} \Po_H^{[N]} - \Po_H^{[N]} \Lo_{-\So_H^{[N]} H_i}  + \Ord(\alpha^{N+1}).\label{phdot1}
\end{equation}

Combining this with  the Lie transform, 
 we can conclude that   the transformation
with generator $G^{[N]}=- \So_H^{[N]} H_i$ normalizes Hamiltonian to the order $N$:
\begin{multline*}{\widetilde H}^{[N]} = \iexpD{\alpha\Lo_{- \So_H^{[N]} H_i}}H + \Ord(\alpha^{N+1})=  \iexpD{\alpha\Lo_{- \So_H^{[N]} H_i}} \Po_H^{[N]} H + \Ord(\alpha^{N+1}) \\ 
= \Po_{H_0}\iexpD{\alpha\Lo_{- \So_H^{[N]} H_i}}  H + \Ord(\alpha^{N+1}) = \Po_{H_0}{\widetilde H}^{[N]} + \Ord(\alpha^{N+1}). 
\end{multline*}
This establishes formally asymptotic character of  series. 

The small denominators also affect  the unperturbed operators $\Po_{H_0}$ and $\So_{H_0}$ defined by \eref{AAdefs}. 
This classical problem was first encountered in celestial mechanics by H. Bruns in 1884. 
Fundamental works of  Kolmogorov \cite{Kolmogorov53} and Arnold \cite{Arnold63} 
thoroughly investigated the analytical  properties of  
averaging operation and the solution of homological equation. 
For general non-degenerate multidimensional  system these operators are analytic for all frequencies, except for a set of Lebesgue measure zero. 
This is sufficient for our formal constructions.
It is worth mentioning that the problem does not concern  applications in which perturbation is represented by finite Fourier sums.

\subsection{Non-uniqueness of the normalized Hamiltonian.}

Similarly to \eref{Final}, 
 the general form of the truncated 
 generator of  the Poincar\'e-Lindstedt transformation is
\[
G^{[N]} = - \So_H^{[N]} H_i + \Po_H^{[N]} F,
\]
where  $F({\bf x})$  defines the normalization style. 
 Now it is possible to discuss  its effects. 

For the sake of clarity, we will  hereafter use formal ``analytic''  expressions, as in the  previous chapters. 
However, one must always remember that for mathematical correctness, these expressions must be straightforwardly converted into 
 truncated sums  (up to   $\Ord(\alpha^{N+1})$).

Following Koseleff \cite{Koseleff}, consider two Hamiltonians normalized in different styles:
$\widetilde H_1= \iexpD{\alpha\Lo_{G_1}}\,H$ and $\widetilde H_2= \iexpD{\alpha\Lo_{G_2}}\,H$. These quantities are connected by the transformation
\[ \widetilde H_2 = \iexpD{\alpha\Lo_{G_2}} \  \expD{\alpha\Lo_{G_1}}\, \widetilde H_1 = \Us_{21} \widetilde H_1.
\]
Obviously, $\Us_{21}$ is canonical. Let us find its generator.  By derivation, we obtain
\begin{multline*} \frac{\partial}{\partial\alpha}\Us_{21} = \iexpD{\alpha\Lo_{G_2}} \left(\Lo_{G_1}-\Lo_{G_2}\right) \expD{\alpha\Lo_{G_1}}
= \iexpD{\alpha\Lo_{G_2}}\, \Lo_{G_1-G_2}\, \expD{\alpha\Lo_{G_1}}\\
= \iexpD{\alpha\Lo_{G_2}}\, \Lo_{\Po_H F_{21}}\, \expD{\alpha\Lo_{G_1}}= \Lo_{\iexpD{\alpha\Lo_{G_2}}\,\Po_H F_{21}}\, \Us_{21}\ .
\end{multline*}
Due to \eref{Final}, the difference of the generators always has the form $\Po_H F_{21}$ with some function $F_{21}({\bf x})$. 
Because of  
\eref{Symeq},  the generator of  $\Us_{21}$ is secular
\[G_{21}=\iexpD{\alpha\Lo_{G_2}}\, \Po_H F_{21} = \Po_{H_0} \iexpD{\alpha\Lo_{G_2}} \, F_{21}.
\]

We may conclude that
{\em the normalized Hamiltonians are connected by the Lie transform, 
 with the generator
belonging to the algebra of integrals of unperturbed system.}

For non-resonance  system with incommensurable frequencies, all these integrals are commutative. 
As a consequence,  {\em non-resonance normal form is unique and  insensitive to normalization style} \cite{Koseleff}. 

This  is not so in the case of resonance. Because  resonance relations give rise to non-commutative  integrals,  the 
{\em  normalized resonance Hamiltonian  depends on the style}. We can obtain  expression 
 for $\tilde H$ following idea of Vittot \cite{Vittot0}. Consider the derivative  
\begin{multline*}
 \frac{\partial}{\partial\,\alpha}\tilde H =  \left(\frac{\partial}{\partial\,\alpha}  \iexpD{\alpha\Lo_{G}}\right) H + \iexpD{\alpha\Lo_{G}}  \frac{\partial\, }{\partial\,\alpha} H =
\iexpD{\alpha\Lo_{G}} \, (\Lo_H G+H_i) \\
= \iexpD{\alpha\Lo_{-\So_H H_i + \Po_H F}} \, \Po_H \, \left(H_i + \D_H F\right) = \Po_{H_0}\, \iexpD{\alpha\Lo_{-\So H_i + \Po_H F}} \, (H_i + \D_H F).
\end{multline*}
Here we  use \eref{inverseD} and the canonical connection of the projectors \eref{Symeq}. Therefore,
\begin{equation}
\tilde H =  H_0 + 
\Po_{H_0} \int_0^\alpha \iexpD{\epsilon\Lo_{-\So_H H_i + \Po_H F}} (H_i + \D_H F) \,\rmd\epsilon.\label{FinalH}
\end{equation}
This  expression demonstrates  the explicit dependence on $F({\bf x})$. 

\subsection{Gustavson integrals}

It is of interest to find physically
meaningful quantities that  are insensitive to an artificial choice of normalization style $F({\bf x})$.
Consider a system with constant unperturbed frequencies.
In his celebrated article 
Gustavson \cite{Gustavson} constructed the formal integrals for a perturbed  system  originating from 
the { centre of algebra of integrals of unperturbed system}.

More precisely, each  
 resonance relation for unperturbed frequencies 
$(\vecp \omega,\vecp D_k)=0$, $k=1,\ldots ,{\textbf r}$ results in an additional non-commutative integral.
In Birkhoff and ``Action-Angle'' representations,
the centre of corresponding  algebra of integrals 
consist of the ${\bf d}-{\bf r}$ quantities
\[\tilde I_i = \sum_{j=1}^{\bf d} \beta_{ij}\tilde\zeta_j \tilde\eta_j = (\vecp\beta_i, \vecp J),\qquad i=1,{\bf d}-{\bf r}.
\]
Here, ${\vecp \beta}_i$ is a set of ${\bf d}-{\bf r}$ independent vectors orthogonal to all ${\bf r}$ resonance vectors
$\vecp D_k$ \cite{Gustavson}.

These integrals are commutative with all  integrals of the unperturbed system, and therefore 
with normalized Hamiltonian ${\widetilde H}$. In operator notation, for any analytic  $\tilde F({\bf x})$,
\begin{align*}
&\tilde I_i = \Po_{H_0} \tilde I_i, \qquad\qquad i=1,{\bf d}-{\bf r},\\
&[ \tilde I_i, \Po_{H_0} \tilde F ]=0,\\
&[ \tilde I_i,  {\widetilde H} ]=0.
\end{align*}

After the transformation back to the initial variables,  the quantities
\begin{equation} I_i = \expD{\alpha\Lo_{G}}\ \tilde I_i,\qquad i=1,{\bf d}-{\bf r},
\end{equation}
 become formal integrals  of the perturbed system. Moreover, these Gustavson integrals are commutative with all functions in the image of $\Po_{H}$: 
\begin{align*}
&I_i = \expD{\alpha\Lo_{G}} \, \Po_{H_0} \tilde I_i =  \Po_{H} \, \expD{\alpha\Lo_{G}} \,\tilde I_i = \Po_{H} I_i, \qquad i=1,{\bf d}-{\bf r},\\
&[ I_i, \Po_{H}  F ]=\expD{\alpha\Lo_{G}} \, [ \tilde I_i, \Po_{H_0} \iexpD{\alpha\Lo_{G}} \, F ] = 0,\\
&[  I_i,  H ]=\expD{\alpha\Lo_{G}} \, [ \tilde I_i, \tilde H ] = 0.
\end{align*}
Here we have again used the canonical connection of the projectors  \eref{Symeq}.
Due to the above properties,  the derivative of $I_i(\alpha)$ will not depend on $F({\bf x})$:
\[
 \left( \frac{\partial}{\partial\,\alpha}I_i (\alpha) \right) =  \left(\frac{\partial}{\partial\,\alpha}  \expD{\alpha\Lo_{G}}\right) \tilde I_i  =
\Lo_{G}  I_i 
= - \Lo_{\So_H H_i}  I_i + \Lo_{\Po_H F} I_i  = - \Lo_{\So_H H_i}  I_i  .
\]
Therefore, Gustavson integrals $I_i(\alpha)$ are insensitive to normalization style.  
Actually, these quantities diverge \cite{Contopoulos2003}, but are useful  in exploring 
the regions of regular dynamics.

The unperturbed Hamiltonian $H_0$ itself may be chosen as a seed for  $\tilde I_i$. 
First orders of the nontrivial part of the corresponding Gustavson integral are
\small\begin{multline*}
I_G =  \alpha^{-1} \left( H - \expD{{\alpha\Lo_{-\So_H H_i}}} H_0 \right) = \Po H_i - \alpha \left( \So\Lo\Po H_i + \tfrac{1}{2}\Po\Lo\So H_i \right)\\
       + \alpha^2 \left( \So \Lo \So \Lo \Po H_i + \tfrac{1}{2} \So \Lo \Po \Lo \So H_i + \tfrac{1}{3} \Po \Lo \So \Lo \So H_i
 - \tfrac{2}{3} \Po \Lo \So^2 \Lo \Po H_i - \tfrac{1}{3} \Po \Lo \Po \Lo \So^2 H_i \right) 
+\Ord(\alpha^{3}) .
\end{multline*}\normalsize
This series is also applicable to general system with non-constant unperturbed frequencies. See  the demo ``{\tt Gustavson\_integral}''. It is known as the Hori's  formal first integral \cite{Hori}.

\paragraph{Example 2. H\'enon-Heiles system.}  This is a two-dimensional system with  Hamiltonian 
\[H= \tfrac{1}{2} (p_1^2+q_1^2 + p_2^2+q_2^2)+ \alpha (q_1^2 q_2-\tfrac{1}{3} q_2^3).
\]
 In  complex $\zeta, \eta$ variables  \eref{zetaeta} the H\'enon-Heiles Hamiltonian becomes
\begin{multline*}H = \rmi  \zeta_2\eta_2 + \rmi \zeta_1 \eta_1
+\alpha \tfrac{1}{2\sqrt{2}} \left( \zeta_2\eta_2^2 - \tfrac{1}{3} \zeta_2^3 
- \eta_1^2\zeta_2 - 2\zeta_1\eta_1\eta_2 +\zeta_1^2\zeta_2 \right.\\
\left. + \tfrac{1}{3} \rmi \eta_2^3 -  \rmi\zeta_2^2\eta_2
          - \rmi \eta_1^2\eta_2 + 2\rmi \zeta_1\eta_1\zeta_2 + \rmi\zeta_1^2\eta_2 \right).
\end{multline*}
The first orders of the normal form are
\begin{multline*}\widetilde H =       \rmi\, (\zeta_2 \eta_2 + \rmi \zeta_1 \eta_1)
       + \alpha^2   \left( \tfrac{5}{12} \zeta_2^2 \eta_2^2 + \tfrac{7}{12} \eta_1^2 \zeta_2^2 - \tfrac{1}{3}\zeta_1 \eta_1 \zeta_2 \eta_2 
 + \tfrac{7}{12}\zeta_1^2 \eta_2^2 + \tfrac{5}{12}\zeta_1^2 \eta_1^2 \right) \\
       + \rmi \alpha^4 \left( \tfrac{235}{432}  \zeta_2^3 \eta_2^3 - \tfrac{175}{144}  \eta_1^2 \zeta_2^3 \eta_2 
- \tfrac{47}{16} \zeta_1 \eta_1 \zeta_2^2 \eta_2^2 + \tfrac{161}{144} \zeta_1 \eta_1^3 \zeta_2^2 \right.\\
\left.- \tfrac{175}{144} \zeta_1^2 \zeta_2 \eta_2^3 + \tfrac{65}{16} \zeta_1^2\eta_1^2\zeta_2\eta_2 
+ \tfrac{161}{144} \zeta_1^3 \eta_1 \eta_2^2 - \tfrac{101}{432}\zeta_1^3 \eta_1^3 \right) + \Ord(\alpha^{6}).
\end{multline*}\normalsize
Because this is the $1:1$ resonance system, we see here mixed terms, such as $\zeta_1 \eta_2$.  
First  orders of the Gustavson  integral are
\small\begin{multline*}I_G= \alpha^{-2} (H - \expD{{\alpha\Lo_{G}}} H_0 ) =  -
\tfrac{1}{48} \left(5 p_1^4+2 p_1^2 \left(5 p_2^2+5 q_1^2-9  q_2^2\right)+56 p_1 p_2 q_1 q_2 +5 p_2^4\vphantom{{q^2}^2}\right.\\\left.
{}-2 p_2^2 \left(9 q_1^2-5 q_2^2\right)
+5\left(q_1^2+q_2^2\right)^2\right) 
-  \alpha\tfrac{1}{36} \left(-28 p_1^4 q_2+28 p_1^3 p_2
   q_1+p_1^2 q_2 \left(84 p_2^2-27 q_1^2+37 q_2^2\right) \vphantom{{q^2}^2}\right.\\\left.
{}+42 p_1 p_2 q_1 \left(-2
   p_2^2+q_1^2+q_2^2\right)-p_2^2 \left(69 q_1^2
   q_2+5 q_2^3\right)-5 q_2 \left(q_2^2-3 q_1^2\right)
   \left(q_1^2+q_2^2\right)\vphantom{{q^2}^2}\right) + \Ord(\alpha^{2})
\end{multline*}\normalsize
In these orders, the expressions are identical to those of classical  works \cite{Gustavson,Jaffe82}. 
Further orders can be found in the demo ``{\tt Henon\_Heiles}''.
The differences between these  and the Deprit series begin at   fifth order in generators and the eighth  order in the normalized Hamiltonians.

As expected, the difference between generators belongs to the kernel of $1-\Po_H$ and 
the normalized Hamiltonians are connected by the Lie transfrom 
with the secular generator.
The Gustavson integrals are identical up to the highest order that we computed.

It is now necessary to say several words about nonresonant systems. 
The KAM theory \cite{Arnold63,Treschev} showed that nonresonant dynamical systems are stable against small perturbations in the sense that
a majority of its invariant tori will not destruct. 
Works of Elliason, Gallavotti and others  \cite{Eliasson,Gallavotti} demonstrated the convergence of Poincar\'e-Lindstedt series  on nonresonant set and
 the systematic cancellations of small denominators. It is natural to search for nonresonance cancellations in Kato series too. 
Such cancellations, connections to Whittakker's adelphic integrals and corresponding computational algorithms  will be discussed in separate article.

\section{Computational aspects.}

A major difference between this and  classical perturbation algorithms by Deprit \cite{Deprit}, Hori \cite{Hori}, Dragt and Finn \cite{DragtFinn}, etc.
is the  explicit non-recursive formulae.  Traditionally, 
perturbation  computations solve homological equations order by order.
In contrast, we   directly compute  generator $G = - \So_H H_i$ up to  the 
desired order  as a sum of all permutations \eref{explicit}.
This  is reminiscent of diagrammatic expansions in quantum field theory.
Then, Lie transform normalizes the Hamiltonian. We present the details of the explicit ``Square''  algorithm for the generator  in Appendix C.

The explicit  expressions are important from a general mathematical point of view, as   
they  systematize and simplify the perturbation expansion,  at the price of additional secular terms in the generator.
Due to  these  terms,  our approach will be  less effective for the practical computations than the classical 
Deprit algorithm.
However, this  should not be a serious problem for contemporary computers.

We compared  the computational times of the direct normalization $\widetilde H= \iexpD{{\alpha\Lo_{-\So_H H_i}}} H$ using the explicit expression \eref{explicit} 
 with those of the Deprit algorithm \cite{Deprit} and the algorithm of Dragt and Finn \cite{DragtFinn}. The latter method uses the products of canonical transformations 
 \[\widetilde H = \rme^{- \alpha^n \Lo_{ G_{n-1}}} \ldots \rme^{- \alpha \Lo_{ G_0}}  H,\qquad\tilde {\bf{x}} = \rme^{ \alpha \Lo_{ G_0}} \ldots \rme^{\alpha^n \Lo_{ G_{n-1}}}   {\bf{x}}, \]
instead of Lie transforms (Deprit exponents) 
in order to normalize Hamiltonian.

Table \ref{tnorm}. compares the times for computing the normal form 
for the one-dimensional anharmonic  oscillator and the two-dimensional H\'enon-Heiles system
on an Intel Xeon X5675 (3.06 GHz) processor.
\begin{table}
\caption{\label{tnorm}Normalization time (s).}
\begin{indented}
\lineup
\item[]\begin{tabular}{@{}lllllll}
\br
&\centre{3}{Anharmonic (1d)}&\centre{3}{H\'enon-Heiles (2d)}\\
\ns
&\crule{3}&\crule{3}\\
\ns
Order&Deprit&Dragt\&Finn&Explicit&Deprit&Dragt\&Finn&Explicit\\
\mr
\04&\00.02&\00.01&\0\00.04&\0\0\00.05&\0\0\0\00.04&\0\0\00.1\\
\08&\00.12&\00.06&\0\00.3&\0\0\00.5&\0\0\0\01.2&\0\0\01.3\\
16&\01.7&\00.8&\0\06.&\0\019.&\0\0139.&\0\042.\\
24&10.7&\03.5&\049.&\0345.&\02231.&\0558.\\
32&46.1&10.2&235.&3155.&20930.&4332.\\
\br
\end{tabular}
\end{indented}
\end{table}

The table reveals that  our method is indeed slower than Deprit, but is still fast enough to complete practical computations. 
 Dragt and Finn's method is faster for one-dimensional systems and lower orders, but is less effective for computing higher orders.
This is a consequence of the growing number of exponentiations intrinsic to it.
The number of the Poisson bracket evaluation in Dragt and Finn’s method is smaller than that of  in Deprit’s algorithm \cite{Cary},
but the brackets have larger arguments.

All the methods build near-identity normalizing canonical transformations.   Such  transformations   are always  the Hamiltonian flows equivalent to Deprit transformation with different normalization styles $\Po_H F$.
As expected,  all three normal forms differ from each other starting from  the  $8\textsuperscript{th}$ order,
while the series for Gustavson integrals coincide. It is the advantage of explicit formula \eref{explicit}   that we are able to compute Gustavson integrals
without the previous normalization.

It is worth  noting that the above computations were single-threaded. If needed, the sum of all the placements  in \eref{explicit} can be  { parallelized} and
made  scalable for contemporary multi-CPU and Cloud computing.

\section{Beyond the perturbation expansion.}
Because  we now have the explicit expression for the generator of the Poincar\'e-Lindstedt transformation, it is interesting to apply it to determine an exact solution
and construct an explicit expression for $\So_H H_i$.
This is possible only for trivial systems because it requires algebraic integrability for all $\alpha$.
However, the following toy examples nicely illustrate the perturbation expansion:
\paragraph{Example 3. Shift of frequency of the harmonic oscillator.}
\[
H=\frac{1}{2}(p^2+q^2)+\frac{\alpha}{2} q^2.\nonumber 
\]
Here we know the exact solution, $q(t)=\sqrt{\frac{2E}{1+\alpha}}\cos(\sqrt{1+\alpha} t+\phi)$,
where  the constant energy $E=H(p,q)$ and $\phi$ are the functions of the initial point $(p_0,q_0)$. 

 We can  directly calculate the ``exact'' perturbation operators using $(\ref{Adef})$ 
and then substitute $p_0\to p$, $q_0\to q$ as follows: 
\begin{align*}
\Po_H \frac{q^2}{2} &= 
\lim_{\lambda\to +0}\lambda\int_0^{+\infty} \rme^{-\lambda t} \frac{q(t)^2}{2}\,\rmd t =
\frac{H}{2(1+\alpha)},\\
\So_H \frac{q^2}{2} &=- \lim_{\lambda\to{+0}}\int_0^\infty \rme^{-\lambda t} \left(\frac{q(t)^2}{2} - \Po_H H_i \right)  \,\rmd t = -\frac{pq}{4(1+\alpha)}.
\end{align*}
 Indeed,
\[\Lo_H (-{pq \over 4(1+\alpha)}) = \frac{-p^2 + q^2 (1 + \alpha)}{4 (1 + \alpha)} =  \frac{q^2}{2} - \Po_H \frac{q^2}{2} .\]
Therefore, the exact generator is $G=  \frac{pq}{4(1+\alpha)}$. 
The normalizing canonical  transformation (``Deprit exponent'') is determined by the equations
\begin{eqnarray}
\frac{\partial p}{\partial\,\alpha} &=& \Lo_G p =-{p\over4(1+\alpha)},\nonumber\\
\frac{\partial q}{\partial\,\alpha} &=& \Lo_G q ={q\over4(1+\alpha)},\nonumber
\end{eqnarray}
with the straightforward solution 
\[
\tilde p ={1\over\sqrt[4]{1+\alpha}}p,\qquad
\tilde q =\sqrt[4]{1+\alpha}q.
\]  The  normalized Hamiltonian is $\widetilde H =\frac{\sqrt{1+\alpha}}{2}({\tilde p}^2 + {\tilde q}^2)$. 
The power series for these exact generator and Hamiltonian coincide with the standard perturbation expansion.

\paragraph{Example 4. The Duffing equation. Part II.}
\[ H=\frac{1}{2}(p^2+q^2)+\frac{\alpha}{4} q^4.
\]
This system allows for an exact solution using Jacobi elliptic functions:
\begin{align*} q(t) &= A \cn(\omega t + \psi, k^2),\\
\omega^2&=\sqrt{1+4\alpha E},\quad
A^2=\frac{1}{\alpha}(\omega^2-1),\quad
k^2=\frac{1}{2}(1-\frac{1}{\omega^2}).
\end{align*}
Here, the energy $E=H(p_0,q_0)$ and ``pseudo-phase'' $\psi$ are determined  by the initial condition
\begin{equation}
\frac{\sn(\psi,k^2)}{\cn(\psi,k^2)}\dn(\psi,k^2) = - \frac{p_0}{q_0\omega}.\label{psi}
\end{equation}

In order to find the exact integrating operator $\So_H  H_i$, 
we may expand the perturbation $H_i=\frac{1}{4}q^4$ into a Fourier series using  \cite{Kiper} 
\begin{equation*}
\cn^4(u,k^2) = {const}
{}+\frac{4\pi^2}{3 k^4\eK^2} \sum_{n=1}^\infty \frac{n Q^n}{1-Q^{2n}} \left((2k^2-1) + \frac{n^2\pi^2}{4K^2}\right) \cos(\frac{n\pi}{\eK}u).
\end{equation*}
Here, $\eK(k^2)$ is the complete elliptic integral of the first kind, and  
\[Q=\exp\left({-\frac{\pi \eK(k^2)}{\eK(1-k^2)}}\right)\]
 is the ``elliptic nome''.

The normalizing generator will be
\begin{multline}\label{q4e}
 G=-\So_H H_i = \lim_{\lambda\to{+0}}\int_0^\infty \rme^{-\lambda t} (\frac{{q(t)}^4}{4} - \Po_H H_i)  \,\rmd t\\
= - \frac{\pi A^4}{3 k^4\omega\eK} \sum_{n=1}^\infty \frac{Q^n}{1-Q^{2n}} \left( \frac{n^2\pi^2}{4K^2}-\frac{1}{\omega^2}\right) \sin(\frac{n\pi}{\eK}\psi).
\end{multline}

For comparison with standard perturbation expansion, 
we  introduce the $2\pi$-periodic ``pseudo-phase'' $\theta=\frac{\pi}{2\eK}\psi$.  This angle can be obtained from  $(\ref{psi})$ as a power series
(again, we substituted $p_0\to p$, $q_0\to q$):
\[
\tan(\theta) = -\frac{p}{q} + \frac{p(3 p^2 + 5 q^2)}{4q(p^2+q^2)} k^2 +\Ord(k^4).
\]

Now we can expand the exact formula $(\ref{q4e})$ into a power series  with the help of Wolfram {\em Mathematica}\textregistered :  
\[
 G=\frac{p q}{32}   (3 p^2  + 5 q^2) -\alpha\frac{p q}{384}   (39 p^4 + 104 p^2 q^2 + 57 q^4)+\Ord(\alpha^2).\\
\]
This  coincides with  perturbation series $(\ref{q4})$. 
Surprisingly, the expansion consumes much more computer time than do perturbative computations.

\section{Summary.}

Being   inspired by the deep parallelism between quantum and classical  perturbation theories, 
we have applied Kato resolvent perturbation expansion to classical mechanics. 
Invariant definitions of averaging and integrating   operators  and  canonical identities uncovered  the regular pattern of perturbation series. 
This pattern was explained using the relation of   perturbation  operators to the Laurent coefficients of Liouville operator resolvent.

 The Kato series for perturbed resolvent and the resolvent canonical identity systematize the perturbation expansion and 
lead to new explicit expression  for 
the  Deprit  generator of Poincare-Lindstedt transformation in any order: 
\begin{displaymath}
G = - \So_H H_i.
\end{displaymath}
Here, the integrating operator  $\So_H$  is the partial pseudo-inverse of the perturbed Liouville operator.  
We have used non-perturbative examples to illustrate this formula.

After extending the  formalism  to multidimensional systems,   
 we  described ambiguities of generator and normalized Hamiltonian. Interestingly, Gustavson integrals turn out to be   insensitive to normalization style. 

All our discussion has remained at a formal level. A comparison of computational times for   this approach 
  and for classic Deprit and Dragt\&Finn algorithms demonstrated that our series is reasonably efficient even for  high orders of perturbation expansion.

\ack
The author gratefully acknowledges  Professor S.~V.\ Klimenko  for reviving my interest in this work and
 RDTEX Technical Support Centre Director S.~P.\ Misiura for encouraging and supporting of my investigations.
 
\appendix
\section{Properties of perturbed operators.}

Let us  explore the properties of truncated perturbed operators $\Po_H^{[N]}$, $\So_H^{[N]}$ and $\D_H^{[N]}$:
\begin{align*}
\Po_H^{[N]}&=\sum_{n=0}^N{(-1)}^{n+1}\alpha^n\left(\sum_{
\sum{p_i}={\bf n}\atop p_i\ge 0}\Sr_{H_0}^{(p_{n+1})}\underbrace{\Lo_{H_i}\Sr_{H_0}^{(p_n)}\,\ldots\,\Sr_{H_0}^{(p_2)}
\Lo_{H_i}}_{n\ {\rm times}}\Sr_{H_0}^{(p_1)}\right),\\
\So_H^{[N]}&=\sum_{n=0}^N{(-1)}^{n\hphantom{+1}}\alpha^n\left(\sum_{
\sum{p_i}={\bf n+1}\atop p_i\ge 0}\Sr_{H_0}^{(p_{n+1})}\underbrace{\Lo_{H_i}\Sr_{H_0}^{(p_n)}\,\ldots\,\Sr_{H_0}^{(p_2)}
\Lo_{H_i}}_{n\ {\rm times}}\Sr_{H_0}^{(p_1)}\right),\\
\D_H^{[N]}&=\sum_{n=1}^N{(-1)}^{n+1}\alpha^n\left(\sum_{
\sum{p_i}={\bf n-1}\atop p_i\ge 0}\Sr_{H_0}^{(p_{n+1})}\underbrace{\Lo_{H_i}\Sr_{H_0}^{(p_n)}\,\ldots\,\Sr_{H_0}^{(p_2)}
\Lo_{H_i}}_{n\ {\rm times}}\Sr_{H_0}^{(p_1)}\right).
\end{align*}
 Here,  ${[N]}$ is the index. We will assume that $N\ge2$. It is useful to introduce operators
\begin{equation}
\Z_n^m = \begin{cases} {(-1)}^{n+1}\sum\limits_{
p_1+\,\ldots\,+p_{n+1}=m\atop p_i\ge 0}\Sr_{H_0}^{(p_{n+1})}
\underbrace{\Lo_{H_i}\Sr_{H_0}^{(p_n)}\,\ldots\,\Sr_{H_0}^{(p_2)}
\Lo_{H_i}}_{n\, {\rm times}}\Sr_{H_0}^{(p_1)},&\mbox{if }m\ge 0,\\
0,&\mbox{if }m<0.\end{cases}
\label{Zrec}
\end{equation}
The summation runs over all possible placements  of $m$  operators $\So_{H_0}$  in $n+1$ sets. 

It is important that $\Z$ operators  can be computed recursively. For  
any $ k, 0\le k \le n-1$
\begin{equation}
\Z_n^m = \sum_{i=0}^m \Z_{n-k-1}^{m-i} \Lo_{H_i} \Z_k^i.\label{Zrecursive}
\end{equation}
The simplest such operators are $\Z_0^0=\Po$, and  $\Z_0^m=-\So^m$ for $m>0$. Therefore, 
\begin{align*}
\Z_0^m \Lo_{H_0} &=\Lo_{H_0} \Z_0^m = \Z_0^{m-1} - \delta_1^m,\\
\Z_0^m H_0\  &= H_0\,\delta_0^m.
\end{align*}
Here $\delta_1^m$ is the Kronecker delta. For
 $n\ge1$  
we have:
\begin{align*}
\Z_n^m \Lo_{H_0} &= \sum_{i=0}^m \Z_{n-1}^{m-i} \Lo_{H_i} \Z_{0}^{i} \Lo_{H_0}  = 
\sum_{i=1}^m \Z_{n-1}^{m-i} \Lo_{H_i} \Z_{0}^{i-1} - \Z_{n-1}^{m-1} \Lo_{H_i} 
= \Z_{n}^{m-1} - \Z_{n-1}^{m-1} \Lo_{H_i},\\
 \Lo_{H_0} \Z_n^m &= \sum_{i=0}^m \Lo_{H_0} \Z_{0}^{i} \Lo_{H_i} \Z_{n-1}^{m-i}  = 
\sum_{i=1}^m \Z_{0}^{i-1}  \Lo_{H_i}  \Z_{n-1}^{m-i} - \Lo_{H_i} \Z_{n-1}^{m-1} 
= \Z_{n}^{m-1} - \Lo_{H_i} \Z_{n-1}^{m-1}.\\
\Z_n^m \; H_0 &= \sum_{i=0}^m \Z_{n-1}^{m-i} \Lo_{H_i} \Z_{0}^{i} H_0 = \Z_{n-1}^{m} \Lo_{H_i}  H_0 = -   \Z_{n-1}^{m} \Lo_{H_0}  H_i = 
- \Z_{n-1}^{m-1} H_i + \delta_1^n\delta_1^m H_i. 
\end{align*}

Next, we will use  finite double-indexed operatorial sums
\begin{equation}
\Sr_{H}^{[N](k)} = -\sum_{n=0}^N \alpha^n \Z_n^{n+k},
\end{equation}
having the following asymptotic behavior
\begin{equation*}
\Sr_{H}^{[N](k)} =\begin{cases}\So^k+\Ord(\alpha)&\mbox{if }k\ge 1,\\
-\Po+\Ord(\alpha)&\mbox{if }k=0,\\
\Ord(\alpha^{k})&\mbox{if }k<0.
\end{cases}
\end{equation*} 
These sums unify the expressions for
\[ \Po_H^{[N]}=-\Sr_{H}^{[N](0)},\qquad  \So_H^{[N]}=\Sr_{H}^{[N](1)},
\qquad  \D_H^{[N]}=-\Sr_{H}^{[N](-1)}.
\]

Their actions on   the perturbed Hamiltonian $H$ are as follows:
\small\begin{align}
\Sr_{H}^{[N](k)} \Lo_{H} & 
=-\Z_0^{k}\Lo_{H_0}-\sum_{n=1}^N \alpha^n \left( \Z_n^{n+k}\Lo_{H_0}+\Z_{n-1}^{n+k-1}\Lo_{H_i}\right)\!+ \Ord(\alpha^{N+1})\nonumber\\
&\qquad\,
=\delta_1^k+\Sr_{H}^{[N](k-1)}\! + \Ord(\alpha^{N+1}),\nonumber\\
 \Lo_{H} \Sr_{H}^{[N](k)} & =\ldots=\delta_1^k+\Sr_{H}^{[N](k-1)}\! + \Ord(\alpha^{N+1}),\label{A4}\\
\Sr_{H}^{[N](k)}\, H & 
= - \Z_0^k  H_0 - \sum_{n=1}^N \alpha^n \left( \Z_n^{n+k} H_0 +\Z_{n-1}^{n+k-1} H_i\right)\!+ \Ord(\alpha^{N+1})
=-\delta^k_0 H + \Ord(\alpha^{N+1}).\nonumber
\end{align}\normalsize
Consider now the resolvent series  truncated at $N\textsuperscript{th}$ order in $\alpha$ and $M\textsuperscript{th}$ order in $z$:
\begin{equation*}
\R_H^{[N,M]}(z)=\sum_{k=-N}^{M+1} \Sr_{H}^{[N](k)}z^{k-1}=-\sum_{k=-N-1}^{M}\!\! z^k\left(\sum_{n=0}^N \alpha^n \Z_n^{n+k+1}\right).
\end{equation*}
$\R_H^{[N,M]}(z)$  is an explicitly meromorphic function of complex $z$. We will use it like a  generating function in probability theory. Powers of $z$ will serve as  placeholders  for simultaneous transformations of  expressions.
Due to the  identities listed above, 
\begin{align*}
\R_H^{[N,M]}(z)\,\Lo_H&=\sum_{k=-N}^{M+1} z^{k-1} \Sr_{H}^{[N](k)}\Lo_H 
=1+z\,\R_H^{[N,M]}(z)-z^{M+1}\Sr_{H}^{[N](M+1)}+ \Ord(\alpha^{N+1}),\\
\Lo_H\,\R_H^{[N,M]}(z)&=1+z\,\R_H^{[N,M]}(z)-z^{M+1}\Sr_{H}^{[N](M+1)}+ \Ord(\alpha^{N+1}).
\end{align*}
Here, the term with $z^{M+1}$ represents the edge effect of truncation. Therefore the following expressions will approximate the identity  operator 
\begin{align*}
&\R_{H}^{[N,M]}(z)(\Lo_H - z) + z^{M+1} \Sr_{H}^{[N](M+1)} =  1 + \Ord(\alpha^{N+1}),\label{almostid}\\
&(\Lo_H - z) \R_{H}^{[N,M]}(z) + z^{M+1} \Sr_{H}^{[N](M+1)} =  1 + \Ord(\alpha^{N+1}).\nonumber
\end{align*}
Corresponding approximate  Hilbert  identity is as follows:
\begin{multline*}
\R_{H}^{[N,M]}(z_1) -\R_{H}^{[N,M]}(z_2) = (z_1-z_2) \R_{H}^{[N,M]}(z_1)\R_{H}^{[N,M]}(z_2) \\- z_1^{M+1}  \Sr_{H}^{[N](M+1)} \R_{H}^{[N,M]}(z_2)+ z_2^{M+1}  \R_{H}^{[N,M]}(z_1)\Sr_{H}^{[N](M+1)}+ \Ord(\alpha^{N+1}).
\end{multline*}

Similarly to \eref{ResCoeff} we can obtain the expression for  the product of the two coefficients 
\begin{multline*}
\Sr_H^{[N](m)}\,\Sr_H^{[N](n)} = \left(\frac{1}{2\pi \rmi}\right)^2
\oint_{|z_1|=\epsilon_1}\oint_{|z_2|=\epsilon_2}
\R_H^{[N,M]}(z_1)\R_H^{[N,M]}(z_2)\,z_1^{-m}\,z_2^{-n}\,\rmd z_1\,\rmd z_2\\
=(\eta_m + \eta_n -1) \Sr_H^{[N](m+n)}
-\eta_{m-M-1} \Sr_{H}^{[N](M+1)} \Sr_H^{[N](m+n-M-1)}\\
+ (1-\eta_{n-M-1})  \Sr_H^{[N](m+n-M-1)}  \Sr_{H}^{[N](M+1)}+ \Ord(\alpha^{N+1})
\end{multline*}\normalsize
For any $m$ and $n$ we can choose the truncation order $M$, so that $M\ge m+n+N$ and $\Sr_H^{[N](m+n-M-1)}=\Ord(\alpha^{N+1})$.
Therefore, 
\begin{equation}
\Sr_H^{[N](m)}\,\Sr_H^{[N](n)} =(\eta_m + \eta_n -1) \Sr_H^{[N](m+n)}+ \Ord(\alpha^{N+1})
\end{equation}
From \eref{A4} and the coefficients of powers of $z$ in the above identity, it follows that:
\begin{itemize}
\item  The finite sum $\Po_H^{[N]}$ is  the approximate projector $\Po_H^{[N]} \Po_H^{[N]}  = \Po_H^{[N]} + \Ord(\alpha^{N+1})$.
\item The sums $\Po_H^{[N]}$, $\So_H^{[N]}$ and $\D_H^{[N]}$ obey approximately  the same identities as the corresponding perturbative operators:
\small\begin{align*}
&\Po_H^{[N]} \So_H^{[N]}  = 0+  \Ord(\alpha^{N+1}), &\So_H^{[N]} \Po_H^{[N]}& = 0+  \Ord(\alpha^{N+1}),  \\
&\D_H^{[N]} \So_H^{[N]}   =  0+  \Ord(\alpha^{N+1}), &\So_H^{[N]} \D_H^{[N]} & =0+  \Ord(\alpha^{N+1}), \\
&\Po_H^{[N]} \D_H^{[N]}   =  \D_H^{[N]}+  \Ord(\alpha^{N+1}), &\D_H^{[N]} \Po_H^{[N]} & =  \D_H^{[N]}+  \Ord(\alpha^{N+1}).
\end{align*}\normalsize
\item These operators act on the perturbed Hamiltonian $H=H_0+\alpha H_i$ as follows:
\small\begin{align*}
&\Po_H^{[N]} H                 = H +  \Ord(\alpha^{N+1}), \\
&\So_H^{[N]} H = 0 +  \Ord(\alpha^{N+1}),
&\D_H^{[N]} H          &= 0 +  \Ord(\alpha^{N+1}),  
\end{align*}\normalsize
\item Interactions of Liouville operator with these sums are as follows: 
\small
\begin{align*}
&\Lo_H \Po_H^{[N]}          = \D_H^{[N]} +  \Ord(\alpha^{N+1}), &\Po_H^{[N]} \Lo_H &= \D_H^{[N]} +  \Ord(\alpha^{N+1}), \\
&\So_H^{[N]}\Lo_H           = 1-\Po_H^{[N]}+  \Ord(\alpha^{N+1}), &\Lo_H \So_H^{[N]} &= 1-\Po_H^{[N]}+  \Ord(\alpha^{N+1}), \\
&\Lo_H \D_H^{[N]}           = \left(\D_H^{[N]}\right)^2+  \Ord(\alpha^{N+1}), &\D_H^{[N]} \Lo_H &= \left(\D_H^{[N]}\right)^2+  \Ord(\alpha^{N+1}).
\end{align*}\normalsize
\item Powers of $\So_H^{[N]}$ and $\D_H^{[N]}$ operators can be expressed  as follows: 
\begin{align*}
&\left(\So_H^{[N]}\right)^k  =\phantom{-}\Sr_{H}^{[N](k)\phantom{{}-{}}\,} +  \Ord(\alpha^{N+1})= - \sum_{n=0}^N \alpha^n\Z_n^{n+k}+  \Ord(\alpha^{N+1}),\\
&\left(\D_H^{[N]}\right)^k    =- \Sr_{H}^{[N](-k)} +  \Ord(\alpha^{N+1}) = \phantom{-}\sum_{n=k}^N \alpha^n\Z_n^{n-k}+  \Ord(\alpha^{N+1}).
\end{align*}
\end{itemize}
See also the demo ``{\tt perturbed\_operators}'' in  the supplementary data files.

\section{Canonical connection of the projectors.}
Now we will use the truncated resolvent 
\begin{equation*}
\R_H^{[N]}(z)=-\frac{1}{z}\Po_H^{[N]} + \sum_{k=0}^{2N} z^k \left(\So_H^{[N]}\right)^{k+1}
\!\!\!-\sum_{k=2}^{N+1} z^{-k} \left(\D_H^{[N]}\right)^{k-1}\!=-\!\!\!\!\sum_{k=-N-1}^{2N}\!\! z^k\left(\sum_{n=0}^N \alpha^n \Z_n^{n+k+1}\right),
\end{equation*}
 to prove the canonical connection of the perturbed and unperturbed  projectors.
In this section we  set the order of truncation in $z$ to $2N$, as it was large enough  to avoid edge effects originating from the expression
\begin{equation}
\R_{H}^{[N]}(z)(\Lo_H - z) + z^{2N+1} \left(\So_H^{[N]}\right)^{2N+1} =  1 + \Ord(\alpha^{N+1}).\label{almostid1}
\end{equation}

In particular, because 
$\left(\D_H^{[N]}\right)^m =  \Ord(\alpha^m)$ and $\So_H^{[N]}\,\D_H^{[N]} =  \Ord(\alpha^{N+1})$, then 
for all integer $1\le k\le N$ we have:
\begin{equation}
\res_{z=0} \left(z^{2N+1}\left(\R_H^{[N]}(z)\right)^k\right)=\Ord(\alpha^{N+1}).\label{res2k}
\end{equation}
Also, for any integer $k\ge 2$:
\[\res_{z=0} \left(\left(\R_H^{[N]}(z)\right)^k\right)=\Ord(\alpha^{N+1}).\]

As expected, the derivative of projector  can be related to the   residue:
\small\begin{multline*} \res_{z=0} \left(\R_H^{[N]}(z) \Lo_{H_i}  \R_H^{[N]}(z)\right) = 
\sum_{k=-N-1}^N \sum_{n=0}^N \sum_{m=0}^N \alpha^{n+m} \Z_n^{n+k+1}  \Lo_{H_i} \Z_m^{m-k} \\
=\sum_{n=0}^N \sum_{m=0}^N \alpha^{n+m} \Z_{n+m+1}^{n+m+1} =\sum_{k=0}^N (k+1) \,\alpha^{k} \, \Z_{k+1}^{k+1} + \Ord(\alpha^{N+1})=
 \left( \frac{\partial}{\partial\alpha}{\Po_H}\right)^{[N]} + \Ord(\alpha^{N+1}).
\end{multline*}\normalsize
There was no loss of accuracy, since we truncated the series after the differentiation. 

Now we can construct the finite analogue of canonical resolvent identity \eref{SimpID2}.
The application of   almost identity operator  \eref{almostid1}
to Poisson bracket $[\R_{H}^{[N]}(z_2)H_i,\R_{H}^{[N]}(z)G]$ 
 results in 
\small\begin{multline}
 \R_H^{[N]}(z) \Lo_{H_i}  \R_H^{[N]}(z) = \Lo_{\R_H^{[N]}(z_2) H_i} \R_H^{[N]}(z) -  \R_H^{[N]}(z) \Lo_{\R_H^{[N]}(z_2) H_i} -
z_2 \R_H^{[N]}(z) \Lo_{\R_H^{[N]}(z_2) H_i}  \R_H^{[N]}(z)\\
 - z^{2N+1}\left( \left(\So_H^{[N]}\right)^{2N+1}   \Lo_{\R_H^{[N]}(z_2) H_i} \R_H^{[N]}\!(z) 
 -\R_H^{[N]}\!(z) \Lo_{\R_H^{[N]}(z_2) H_i} \left(\So_H^{[N]}\right)^{2N+1}\right) \\
+ z_2^{2N+1}\R_H^{[N]}\!(z) \Lo_{\left(\So_H^{[N]}\right)^{2N+1}\! H_i} \R_H^{[N]}\!(z) + \Ord(\alpha^{N+1}).\label{finiteres}
\end{multline}\normalsize
We are interested in the residue of the above formula at $z=0$:
\small\begin{multline*}
 \left( \frac{\partial}{\partial\alpha}{\Po_H}\right)^{[N]}= 
\Po_H^{[N]} \Lo_{\R_H^{[N]}(z_2) H_i}  - \Lo_{\R_H^{[N]}(z_2) H_i} \Po_H^{[N]} 
- z_2 \res_{z=0}  \left(\R_H^{[N]}(z) \Lo_{\R_H^{[N]}(z_2) H_i}  \R_H^{[N]}(z)\right)\\
+ z_2^{2N+1}\res_{z=0}  \left(\R_H^{[N]}\!(z) \Lo_{\left(\So_H^{[N]}\right)^{2N+1}\! H_i} \R_H^{[N]}\!(z)\right) + \Ord(\alpha^{N+1}).
\end{multline*}\normalsize
Its coefficient of $z_2^0$ is
\small\begin{multline}
 \left( \frac{\partial}{\partial\alpha}{\Po_H}\right)^{[N]} =     \Po_H^{[N]} \Lo_{\So_H^{[N]} H_i} - \Lo_{\So_H^{[N]}  H_i} \Po_H^{[N]} 
+ \res_{z=0}  \left(\R_H^{[N]}(z) \Lo_{\Po_H  H_i}  \R_H^{[N]}(z)\right) + \Ord(\alpha^{N+1}).\label{lastres}
\end{multline}\normalsize
We will show that the last residue    is $\Ord(\alpha^{N+1})$ also. The coefficient of  $z_2^{-1}$ in \eref{finiteres}  is
\small\begin{multline*}
0=\R_H^{[N]}(z) \Lo_{\Po_H^{[N]} H_i} - \Lo_{\Po_H^{[N]} H_i} \R_H^{[N]}(z) +  \R_H^{[N]}(z) \Lo_{\D_H^{[N]} H_i}  \R_H^{[N]}(z)\\
 + z^{2N+1}\left( \left(\So_H^{[N]}\right)^{2N+1}    \Lo_{\Po_H^{[N]} H_i} \R_H^{[N]}(z)  - \R_H^{[N]}(z) \Lo_{\Po_H^{[N]} H_i} \left(\So_H^{[N]}\right)^{2N+1}\right) + \Ord(\alpha^{N+1}).
\end{multline*}\normalsize
Therefore,
\small\begin{multline*}
\res_{z=0}  \left(\R_H^{[N]}(z) \Lo_{\Po_H^{[N]}  H_i}  \R_H^{[N]}(z)\right) = \res_{z=0}  \left(\Lo_{\Po_H^{[N]} H_i} \left(\R_H^{[N]}(z)\right)^2 \right)
-  \res_{z=0}  \left(\R_H^{[N]}(z) \Lo_{\D_H^{[N]} H_i}  \left(\R_H^{[N]}(z)\right)^2\right)\\
 - \res_{z=0}  \left(z^{2N+1}\left( {\left(\So_H^{[N]}\right)}^{2N+1}    \Lo_{\Po_H^{[N]} H_i} \R_H^{[N]}(z)  - \R_H^{[N]}(z) \Lo_{\Po_H^{[N]} H_i} \left(\So_H^{[N]}\right)^{2N+1}\right)\R_H^{[N]}(z)\right)  + \Ord(\alpha^{N+1}) \\
=- \res_{z=0}  \left(\R_H^{[N]}(z) \Lo_{\D_H^{[N]}  H_i}  \left(\R_H^{[N]}(z)\right)^2\right)
 + \Ord(\alpha^{N+1}).
\end{multline*}\normalsize
We can continue this process using  coefficients of  $z_2^{-2}$, $z_2^{-3}$, \ldots in \eref{finiteres}
\begin{multline*}\res_{z=0}  \left(\R_H^{[N]}(z) \Lo_{\D_H^{[N]}  H_i}  \left(\R_H^{[N]}(z)\right)^2\right) =
 - \res_{z=0}  \left(\R_H^{[N]}(z) \Lo_{\left(\D_H^{[N]}\right)^2  H_i}  \left(\R_H^{[N]}(z)\right)^3\right) + \Ord(\alpha^{N+1}) \\
\ldots 
=(-1)^{N+1}\res_{z=0}  \left(\R_H^{[N]}(z) \Lo_{{\left(\D_H^{[N]}\right)}^N  H_i}  \left(\R_H^{[N]}(z)\right)^{N+1}\right) + \Ord(\alpha^{N+1}) = \Ord(\alpha^{N+1}).
\end{multline*}
Inserting this into \eref{lastres} concludes the proof that the projector canonically transforms  
\[
 \left( \frac{\partial}{\partial\alpha}{\Po_H}\right)^{[N]} =    \Lo_{-\So_H^{[N]}  H_i} \Po_H^{[N]} - \Po_H^{[N]} \Lo_{-\So_H^{[N]} H_i}  + \Ord(\alpha^{N+1}),
\]
with the generator $G^{[N]} = - \So_H^{[N]}  H_i$ in any perturbation order.

\section{The square algorithm.}
Relation \eref{Zrecursive} leads to an efficient  computational algorithm for the generator $G^{[N]} = - \So_H^{[N]}  H_i$. 
Consider a square table:
\[
\begin{array}{lllll}
 F^{0}_{0}({\bf x})       & F^{1}_{0}({\bf x}) &F^{2}_{0}({\bf x})&\ldots&F^{N+1}_{0}({\bf x})\\
 F^{0}_{1}({\bf x})       &F^{1}_{1} ({\bf x})       &F^{2}_{1}({\bf x})&\ldots&F^{N+1}_{1}({\bf x})\\
&\ldots&&&\\
F^{0}_{N-1} ({\bf x})      &F^{1}_{N-1} ({\bf x})      &F^{2}_{N-1}({\bf x})&\ldots&F^{N+1}_{N-1}({\bf x})\\
&&&&F^{N+1}_{N}({\bf x})
\end{array}
\]
Here the first row is as follows:
\[F_0^0({\bf x}) = \Z_0^0 H_i=\Po H_i, \quad F_0^m({\bf x}) =\Z_0^m H_i=-\So^m H_i,\]
and  each next row is generated from the previous one according to the rule
\[ F_{n+1}^m ({\bf x}) = \sum_{i=0}^m \Z_{0}^{m-i} \Lo_{H_i} F_{n}^i({\bf x}). \]
 The normalizing generator is given by $G^{[N]}= \sum_{n=0}^N \alpha^n F_n^{n+1}  ({\bf x}) $. 
Obvious computational optimization is to store the quantities $\Lo_{H_i} F_{n}^i({\bf x})$.
\section*{References}
\bibliographystyle{iopart-num}
\bibliography{kato_expansion}

\end{document}